%%%%%%%%%%%%%%%%%%%%%%%%%%%%%%%%%%%%%%%%%%%%%%%%%
%
% This  PLAIN   TeX   file is called ChanilloKiesslingDMJ.tex
%
%   original May 23, 1999; revised January 13, 2000
%
%%%%%%%%%%%%%%%%%%%%%%%%%%%%%%%%%%%%%%%%%%%%%%%%%%
%
% This file already includes  ``auto-numbering.tex''
%
%%%%%%%%%%%%%%%%%%%%%%%%%%%%%%%%%%%%%%%%%%%%%%%%%%%

\newcount\chno  
\newcount\equno 
\newcount\refno 

\font\chhdsize=cmbx12 at 14.4pt

% Bibliography.  Note:  the bibliography information should come at the
% beginning of the paper, though it will not appear where it is typed.
% Usage:
%   Start the bibliography with \startbib
%   Put entries in bibliography in the desired order.  Numbering is done
%   automatically.  For example
%     \reflbl\donovan{Donovan, G. C., ``Fractal Wavelets'', Bla, bla, bla.}
%   End the bibliography with \endbib
%   Reference bibliography entries; for example
%     ... as was shown in [\donovan].
%   Place the bibliography where you want it to appear in the output by
%     by invoking \biblio
%
\def\startbib{\def\biblio{\bigskip\medskip
  \noindent{\chhdsize References.}\bgroup\parindent=2em}}
\def\endbib{\edef\biblio{\biblio\egroup}}
\def\reflbl#1#2{\global\advance\refno by 1
  \edef#1{\number\refno}
    \global\edef\biblio{\biblio\par\item{[\number\refno]}#2\par}}

% For usage of \eqlbl, see example in file  ref_num_sample.tex
\def\eqlbl#1{\global\advance\equno by 1
  \global\edef#1{{\number\chno.\number\equno}}
  (\number\chno.\number\equno)}

\def\qed{\hfil\hbox to 0pt{}\ \hbox to 2em{\hss}\ 
         \hbox to 0pt{}\hskip-2em plus 1fill
         \vrule height6pt depth1pt width7pt\par\medskip}
\def\eqed{\hfil\hbox to 0pt{}\ \hbox to 2em{\hss}\ 
          \hbox to 0pt{}\hskip-2em plus 1fill
          \vbox{\hrule height .25pt depth 0pt width 7pt
            \hbox{\vrule height 6.5pt depth 0pt width .25pt
              \hskip 6.5pt\vrule height 6.5pt depth 0pt width .25pt}
            \hrule height .25pt depth 0pt width 7pt}\par\medskip}

% \msimp superimposes two characters.  The width of the result is the
% width of the first argument.

\def\msimp#1#2{#1%
  \hbox to 0pt{\hskip 0pt minus 3fill
    \phantom{$#1$}\hbox to 0pt{\hss$#2$\hss}\phantom{$#1$}%
    \hskip 0pt minus 1fill}}

%%%%%%%%%%%%%%%%%%%%%%%%%%

\def\pr{{\prime}}
\def\ppr{{\prime\prime}}

\def\({\left(}
\def\){\right)}
\def\[{\left[}
\def\]{\right]}
\def\<{\left\langle}
\def\>{\right\rangle}

\def\ol{\overline}
\def\ul{\underline}

\def\wtilde{\widetilde}
\def\what{\widehat}

\def\eps{\epsilon}

\font\tendouble=msbm10 
\font\sevendouble=msbm7  % changed by MK 
\font\fivedouble=msbm5

\newfam\dbfam
\textfont\dbfam=\tendouble \scriptfont\dbfam=\sevendouble
\scriptscriptfont\dbfam=\fivedouble

\mathchardef\dbA="7041 
\mathchardef\dbB="7042 
\mathchardef\dbC="7043 \def\CC{{\fam=\dbfam\dbC}}
\mathchardef\dbD="7044 
\mathchardef\dbE="7045 
\mathchardef\dbF="7046 
\mathchardef\dbG="7047 
\mathchardef\dbH="7048 
\mathchardef\dbI="7049 
\mathchardef\dbJ="704A 
\mathchardef\dbK="704B 
\mathchardef\dbL="704C 
\mathchardef\dbM="704D 
\mathchardef\dbN="704E \def\NN{{\fam=\dbfam\dbN}}
\mathchardef\dbO="704F 
\mathchardef\dbP="7050 
\mathchardef\dbQ="7051 
\mathchardef\dbR="7052 \def\RR{{\fam=\dbfam\dbR}}
\mathchardef\dbS="7053 \def\SS{{\fam=\dbfam\dbS}}
\mathchardef\dbT="7054 
\mathchardef\dbU="7055 
\mathchardef\dbV="7056 
\mathchardef\dbW="7057 
\mathchardef\dbX="7058 
\mathchardef\dbY="7059 
\mathchardef\dbZ="705A

\def\cE{{\cal E}}
\def\cF{{\cal F}}
\def\cH{{\cal H}}
\def\cK{{\cal K}}

\def\cP{{\cal P}}

\def\cR{{\cal R}}
\def\cS{{\cal S}}

\def\dtau{\,{\rm d}\tau}
\def\dmu{\,{\rm d}\mu}
\def\dvarpi{\,{\rm d}\varpi}
\def\dtheta{\,{\rm d}\theta}
\def\deta{\,{\rm d}\eta}

\def\dvarrho{\,{\rm d}\varrho}
\def\dlambda{\,{\rm d}\lambda}
\def\dxi{\,{\rm d}\xi}

\def\dr{\,{\rm d}r}
\def\ds{\,{\rm d}s}

\def\dx{\,{\rm d}x}
\def\dy{\,{\rm d}y}

\def\Gave#1{ {\rm Ave} \left( {#1}\right) }

%%%%%%%%%%%%%%%%%%%%%%%%%%%%%

\magnification= 1200 
\baselineskip= 12pt

%%%%%%%%%%%%%%

\input epsf 

%%%%%%%%%%%%%%

\startbib
\refno=0

%%%%%%%%%% Bibliography

\reflbl\Ahlfors{
Ahlfors, L.V., 
	{\it An extension of Schwarz's Lemma},
		Trans. Amer. Math. Soc. {\bf 43}, pp. 359--364 (1938).}

%\reflbl\Aly{
%Aly, J.J., 
%      	{\it Some properties of a thermodynamic model for the
%	equilibrium of a current-carrying quasineutral plasma,}
%        	Phys. Fluids {\bf B5}, 281--294 (1993).}

\reflbl\Aubin{
Aubin, T., 
	{\it Meilleures constantes dans le th\'eo\-r\`eme d'inclusion
	de Sobolev et un th\'eo\-r\`eme de Fredholm non lin\'eaire
        pour la transformation conforme de la courbure scalaire},
        	J. Functional Anal. {\bf 32}, pp. 148--174 (1979).}

\reflbl\Aviles{
Aviles, P.,
	{\it Conformal complete metrics with prescribed non-negative
	Gaussian curvature in $\RR^2$},
		Invent. Math. {\bf 83}, pp 519--544 (1986).}

%\reflbl\Bandle{
%Bandle, C.,
%	{\it Isoperimetric Inequalities and Applications,}
%		Pitman, Boston, (1980).}

\reflbl\Beckner{
Beckner, W., 
	{\it Sharp Sobolev inequalities on the sphere and the
		Moser-Trudinger inequality}, 
		Ann. Math. {\bf 138}, pp. 213--242 (1993).}

\reflbl\BillingsleyBOOK{
Billingsley, P.,
	{\it Convergence of Probability Measures},
		J. Wiley and Sons, New York (1968).}

\reflbl\Branson{
Branson, T.,
	{\it Group representations arising from Lorentz conformal
	geometry}, J. Functional Anal. {\bf 74}, pp. 199--293 (1987).}

\reflbl\BrezisMerle{ 
Brezis, H., and  Merle, F., 
	{\it Uniform estimates and blow-up behavior of solutions of 
	$-\Delta u = V(x) e^u$ in two dimensions,}
		Commun.\ PDE {\bf 16}, pp. 1223--1253 (1991).}

\reflbl\clmpCMPcan{
Caglioti, E., Lions, P. L., Marchioro, C., and Pulvirenti, M.,
	{\it A special class of stationary flows for two-dimensional
	Euler equations: A statistical mechanics description}, 
		Commun. Math. Phys. {\bf 143},  pp. 501--525 (1992).}

\reflbl\CarlenLossB{
Carlen, E., and Loss, M., 
	{\it Competing symmetries, the logarithmic HLS inequality, 
	and Onofri's inequality on $S^n$},
		Geom. Funct. Anal. {\bf 2} pp. 90--104 (1992).}

\reflbl\ChangYangA{
Chang., S.Y.A., and Yang, P., 
	{\it Prescribing Gaussian curvature on $\SS^2$,}
		Acta Math. {\bf 159}, pp. 215--259 (1987).}

\reflbl\ChangYangB{
Chang., S.Y.A., and Yang, P., 
	{\it Conformal metrics on $\SS^2$},
		J. Diff. Geom. {\bf 27}, pp. 256--296 (1988).}

\reflbl\ChaKieCMP{
Chanillo, S., and Kiessling, M. K.-H.,
	{\it Rotational symmetry of solutions of some nonlinear 
	problems in statistical mechanics and geometry}, 
		Commun. Math. Phys. {\bf 160}, 217--238 (1994). }

\reflbl\ChaKieGAFA{
Chanillo, S., and Kiessling, M. K.-H.,
	{\it Conformally invariant systems of nonlinear 
	PDE of Liouville type},
		Geom. Functional Anal., {\bf 5}, pp. 924--947 (1995).}

\reflbl\ChaLi{
Chanillo, S., and Li, Y. Y., 
       	{\it Continuity of solutions of uniformly 
        elliptic equations in $\RR^2$},
	        Manuscr. Math. {\bf 77}, pp.  415--433 (1992).}

\reflbl\ChenLiA{
Chen, W., and  Li, C., 
	{\it Classification of solutions of some nonlinear 
	elliptic equations,}
		Duke Math.\ J. {\bf 63},  pp. 615--622 (1991).} 

\reflbl\ChenLiB{
Chen, W., and Li, C.,
	{\it Qualitative Properties of solutions to some 
	nonlinear elliptic equations in $\RR^2$,}
		Duke Math.\ J. {\bf 71},  pp. 427--439 (1993).} 

\reflbl\ChengLinMA{
Cheng, K.-S., and Lin, C.-S.,
	{\it On the asymptotic behavior of solutions of the 
	conformal Gaussian curvature equations on $\RR^2$}, 
		Math. Ann. {\bf 308}, pp. 119--139 (1997).}

\reflbl\ChengLinJDE{
Cheng, K.-S., and Lin, C.-S.,
	{\it On the conformal Gaussian curvature equations in $\RR^2$}, 
		J. Diff. Eq. {\bf 146}, pp. 226--250 (1998).}

\reflbl\ChengLinASNSP{
Cheng, K.-S., and Lin, C.-S.,
	{\it Compactness of conformal metrics with positive Gaussian 
	curvature in $\RR^2$}, 
		Ann. Scuola Norm. Sup. Pisa Cl. Sci. {\bf 26}, 
		pp. 31--45 (1998).}

\reflbl\ChengLinNA{
Cheng, K.-S., and Lin, C.-S.,
	{\it Conformal metrics in $\RR^2$ with prescribed
	Gaussian curvature with positive total curvature},
		Nonl. Anal. {\bf 38}, pp. 775-783 (1999).}

\reflbl\ChengNiI{
Cheng, K.-S., and Ni, W.-M.,
	{\it On the structure of the conformal Gaussian 
	curvature equation on $\RR^2$}, 
		Duke Math.\ J. {\bf 62},  pp. 721--737 (1991).} 

\reflbl\ChengNiII{
Cheng, K.-S., and Ni, W.-M.,
	{\it On the structure of the conformal Gaussian 
	curvature equation on $\RR^2$, II}, 
		Math. Ann. {\bf 290},  pp. 671--680 (1991).} 

\reflbl\ChouWan{
Chou, K.S., and Wan, T.Y.H.,
	{\it Asymptotic radial symmetry for
	solutions of $\Delta u+e^u=0$ in a punctured disc,}
		Pac. J. Math. {\bf 163}, pp.269--276 (1994).}

\reflbl\Diaconis{
Diaconis, P.,
	{\it Recent progress on de Finetti's notions of exchangeability},
		Bayesian Stat. {\bf 3}, pp. 111--125 (1988).}

\reflbl\DiaconisFreedman{
Diaconis, P., and Freedman, D., 
	{\it A dozen de Finetti-style results in search of a theory},
		Ann. Inst. H. Poincar\'e {\bf 23}, pp. 397--423 (1987).}

\reflbl\DJLW{
Ding, W., Jost, J., Li, J., and Wang, G., 
	{\it The differential equation $\Delta u = 8\pi -8\pi he^u$
	on a compact Riemann surface}, 
		Univ. Leipzig Preprint, Feb. 1998.}

\reflbl\Dynkin{
Dynkin, E. B., 
		{\it Klassy ekvivalentnyh slu$\check{c}$a$\check{i}$nyh 
	veli$\check{c}$in}, 
        Uspeki Mat. Nauk. {\bf 6}, pp. 125--134 [1953].}

\reflbl\EllisBOOK{
Ellis, R. S.,
		{\it Entropy, large deviations, and statistical mechanics}, 
	Springer-Verlag, New York [1985].}

%\reflbl\Fadeev{
%Fadeev, V. M., Kvartskhava, I. F., and Kamarov, N. N., 
%      	{\it Self-focussing of local plasma currents,} (in Russian), 
%	        Nuclear Fusion {\bf 5}, 202--209 (1965). }

\reflbl\Finetti{
de Finetti, B., 
		{\it Funzione caratteristica di un fenomeno aleatorio}, 
	Atti della R. Accad. Naz. dei Lincei, Ser. 6, Memorie, 
	Classe di Scienze Fisiche, Matematiche e Naturali 4,
        pp. 251--300 [1931].}

\reflbl\GidasNiNirbrg{ 
Gidas, B., Ni, W.-M., and Nirenberg, L.,
	{\it Symmetry and related properties via the maximum principle,}
		Commun.\ Math.\ Phys. {\bf 68}, pp. 209--243 (1979).}

\reflbl\GilbargTrudinger{
Gilbarg, D. and  Trudinger, N.S., 
	{\it Elliptic Partial Differential Equations of Second Order,}
		 Springer Verlag, New York (1983).}

\reflbl\GlimmJaffe{
Glimm, J., and Jaffe, A.,
	{\it Quantum Physics}, $2^{nd}$ ed., 
		 Springer Verlag, New York (1987).}

\reflbl\Han{
Han, Z.C.,
	{\it Prescribing Gaussian curvature on $\SS^2$},
		Duke Math. J. {\bf 61}, pp. 679--703 (1990).}

\reflbl\HardyLittlewoodPolya{
Hardy, G., Littlewood, J.E., and P\'olya, G.,
	{\it Inequalities}, 2.ed., 
		Cambridge Univ. Press, Cambridge [1994].}

\reflbl\HewittSavage{
Hewitt, E., and Savage, L. J., 
		{\it Symmetric measures on Cartesian products}, 
        Trans. Amer. Math. Soc. {\bf 80}, pp. 470--501 [1955].}

\reflbl\KazdanWarner{
Kazdan, J.L., and Warner, F.W.,
	{\it Curvature functions for compact 2-manifolds},
		Ann. Math. {\bf 99}, pp. 14--42 (1974).}

\reflbl\KieCPAM{
Kiessling, M.K.-H., 
	{\it Statistical mechanics of classical particles with 
	logarithmic interactions},
		Commun. Pure Appl. Math. {\bf 46}, pp. 27--56 (1993).}

%\reflbl\KieLMP{
%Kiessling, M.K.-H., 
%	{\it Negative temperature bounds for 2D vorticity compounds},
%		Lett. Math. Phys. {\bf 34}, pp. 49--57 (1995).}

\reflbl\KiePHYa{
Kiessling, M.K.-H., 
	{\it Statistical mechanics approach to some problems in conformal
	geometry},
		Physica A, to appear (2000).}

\reflbl\KieJllPHPL{
Kiessling, M.K.-H., and Lebowitz, J.L.,
	{\it Dissipative stationary plasmas: Kinetic modeling, 
	 Bennett's pinch and  generalizations,}
		The Physics of Plasmas {\bf 1}, pp. 1841--1849 (1994).}

\reflbl\KieJllLMP{
Kiessling, M.K.-H., and Lebowitz, J.L.,
	{\it The microcanonical point vortex ensemble: Beyond equivalence}, 
		Lett. Math. Phys.  {\bf 42}, pp. 43--56 (1997).}
	
\reflbl\KieSpo{
Kiessling, M.K.-H., and Spohn, H.,
	{\it A note on the eigenvalue density of random matrices},
		Comm. Math. Phys. {\bf 199}, pp. 683--695 (1999).}

\reflbl\Liouville{
Liouville, J., 
	{\it Sur l'\'equation aux diff\'erences partielles
        $\partial^2 \log \lambda/\partial u\partial v \pm \lambda/ 2 a^2 =0$},
	        J. de Math. Pures Appl. {\bf 18}, pp. 71--72 (1853).}

\reflbl\LiCM{
Li, C.-M., 
	{\it Monotonicity and symmetry of solutions of fully 
	nonlinear elliptic equations on unbounded domains,}
		Commun.\ PDE {\bf 16}, pp. 585--615 (1991). }

\reflbl\LiYY{
Li, Y.-Y., 
	{\it Prescribing scalar curvature on $\SS^n$ and related
	problems, part I},
		Rutgers Univ. Preprint.}

%\reflbl\Lions{
%Lions, P.-L., 
%	{\it Two geometrical properties of solutions of semilinear problems,}
%		Appl.\ Anal. {\bf 12}, pp. 267--272 (1981).}

\reflbl\Ma{
Ma, Li, 
	{\it Bifurcation in the scalar curvature problem}, 
		Rutgers Univ. Preprint Feb. 14, 1996. }

\reflbl\McOwen{
McOwen, R.C.,
	{\it Conformal metrics in $\RR^2$ with prescribed Gaussian
	curvature and positive total curvature},
		Indiana Univ. Math. J. {\bf 34}, pp. 97--104 (1985).}

\reflbl\MoserA{
Moser, J., 
	{\it A sharp form of an inequality by N. Trudinger},
	        Indiana Univ. Math. J. {\bf 20}, pp. 1077--1092 (1971).}

\reflbl\MoserB{
Moser, J., 
	{\it On a nonlinear problem in differential geometry},
		in {\it Dynamical Systems} (M. Peixoto, Ed.),
		Acad. Press, New York (1973).}

\reflbl\Ni{
Ni, W.M., 
	{\it On the elliptic equation $\Delta u + Ke^{2u} =0$ and
	conformal metrics with prescribed Gaussian curvatures},
		Invent. Math. {\bf 66}, pp. 343--352 (1982).}

%\reflbl\Nirenberg{
%Nirenberg, L.,
%	{\it Topics in nonlinear functional analysis}, 
%		Lecture Notes, Courant Inst. Math. Sci., NYU (1973/74).}

\reflbl\Obata{
Obata,  M., 
	{\it The conjectures on conformal transformations of 
	Riemannian manifolds,}
		J.\ Diff.\ Geom. {\bf 6}, pp. 247--258 (1971).}

\reflbl\Onofri{
Onofri, E., 
	{\it On the positivity of the effective action in a 
	theory of random surfaces}, 
		Commun. Math. Phys. {\bf 86}, pp. 321--326 (1982).}

\reflbl\Osserman{
Osserman, R., 
	{\it On the inequality $\Delta u \geq f(u)$},
		Pac. J. Math. {\bf 7} pp. 1641--1647 (1957).}

\reflbl\Poincare{
Poincar\'e, H., 
       	{\it Les fonctions fuchsiennes et l'\'equation $\Delta u = e^u$},
	        J. Math. {\bf 4}, pp. 137--230 (1898).}

\reflbl\PrajapatTarantello{
Prajapat, J., and Tarantello, G.,
	{\it On a class of elliptic problems in $\RR^2$: symmetry and
	uniqueness results},
	Proc. Roy. Soc. Edinburgh (to appear).}

\reflbl\RobinsonRuelle{
Robinson, D. W., and Ruelle, D., 
	{\it Mean entropy of states in classical statistical mechanics}, 
		Commun. Math. Phys. {\bf 5}, pp. 288--300 [1967].}

\reflbl\Rockafellar{
Rockafellar, R.T.,
		{\it Convex analysis},
	Princeton Univ. Press (1970).}
  
\reflbl\RuelleBOOK{
	Ruelle, D., 
		{\it Statistical Mechanics: Rigorous Results},
	Addison Wesley [1989].}

\reflbl\SaffTotikBOOK{
	Saff, E.B., and Totik, V.,
		{\it Logarithmic potentials with external fields},
	Grundl. d. Math. Wiss. vol. {\bf 316}, Springer, New York (1997).}

\reflbl\Sattinger{
Sattinger, D.H.,
	{\it Conformal metrics in $\RR^2$ with prescribed curvatures},
		Indiana Univ. Math. J. {\bf 22}, pp. 1--4 (1972).}

%\reflbl\SchmidtBurgk{
%Schmid-Burgk, J., 
%       	{\it Finite amplitude density variations 
%        in a self-gravi\-tat\-ing isothermal gas layer,}
%	        Astrophys. J. {\bf 149}, pp. 727--729 (1967).}

%\reflbl\SpohnBOOK{
%Spohn, H., 
%	{\it Large Scale Dynamics of Interacting Particles}, 
%		Texts and Monographs in Physics, Springer, New York (1991).}

%\reflbl\Stuart{
%Stuart, J. T., 
%        {\it On finite amplitude oscillations in laminar mixing layers},
%	        J. Fluid Mech. {\bf 29} pp. 417--439 (1967).}

\reflbl\Tarantello{
Tarantello, G.,
	{\it Vortex condensation of a non-relativistic Chern-Simons theory},
		J. Diff. Eq. {\bf 141}, pp. 295--309 (1997).}

\reflbl\Wittich{
Wittich, H., 
	{\it Ganze L\"osungen der Differentialgleichung $\Delta u=e^u$},
		Math. Z. {\bf 49}, pp. 579--582  (1944).}

%\reflbl\Walker{
%Walker, G. W.,
%	{\it Some problems illustrating the forms of nebulae},
%		Proc. Roy. Soc. {\bf A 91}, pp. 410--420 (1915).}

\endbib

%%%%%%%%%%%%%%%%%%%%%%%%%%

\centerline{\chhdsize SURFACES WITH RADIALLY SYMMETRIC }
\medskip
\centerline{\chhdsize    PRESCRIBED GAUSS CURVATURE}
\bigskip
\bigskip
\bigskip
\centerline{SAGUN CHANILLO and MICHAEL KIESSLING}
\bigskip
\centerline{Department of Mathematics}
\centerline{Rutgers University}
\centerline{Piscataway, NJ 08854}
\bigskip
\bigskip
\bigskip
\bigskip
\bigskip

\noindent
{\bf ABSTRACT}: 
We study conformally flat  surfaces with prescribed 
Gaussian curvature, described by solutions $u$ of the PDE: 
$\Delta u(x)+K(x)\exp(2u(x))=0$, with $K(x)$ the Gauss curvature
function at $x\in\RR^2$. We assume that the integral curvature is 
finite. For radially symmetric $K$ we introduce the notion 
of a least integrally curved surface, and also the notion of when 
such a surface is critical. With respect to these notions we analyze 
the radial symmetry of $u$ for the whole spectrum of possible integral
curvature values. Under a mild integrability condition which rules out 
harmonic non-radial behavior near infinity, we prove that $u$ is
radially symmetric and decreasing in the following categories: 
(1) $K$ is decreasing, $u$ a classical solution, and the integral 
curvature of the surface is above critical; 
(2) $K$ is decreasing, $u$ a classical solution, 
the integral curvature of the surface is critical,
and the surface satisfies an additional integrability
condition which is mildly stronger than finite integral curvature;
(3) $K$ is non-positive. 
In categories 1 and 2, $K$ is allowed to diverge logarithmically or
as power law to $-\infty$ at spatial infinity.
Examples of nonradial solutions which violate one or more of our 
conditions are discussed as well. In particular, for non-positive
and non-negative $K$ that satisfy appropriate integrability
conditions and otherwise are fairly arbitrary, we introduce
probabilistic methods to construct surfaces with finite 
integral curvature and entire harmonic asymptotics at infinity. 
For radial symmetric $K$ these surfaces are examples of broken symmetry. 

\hfill

\bigskip\bigskip

\centerline{To appear in Duke Mathematical Journal}

\centerline{Original May 23, 1999; revised Jan. 13, 2000.}
\vfill

\hrule\medskip
\noindent
$\msimp{\copyright}{c}$ (2000) The authors. 
Reproduction for non-commercial purposes of 
this article in its entirety, by any means, is permitted.

\vfill\eject

\chno=1
\equno=0
\noindent
{\bf I. INTRODUCTION}
\medskip
Let $S_g = (\RR^2,g)$ denote a conformally flat surface over $\RR^2$ 
with metric given by 
$$
{\ds}^2 = g^{ij}\dx_i\dx_j = e^{2u(x)}\left({\dx_1}^2 + {\dx_2}^2\right)\, ,
\eqno\eqlbl\metric
$$
where $u$ is a real-valued function of the isothermal coordinates
$x = (x_1, x_2)\in \RR^2$. If $u$ is given, the Gauss curvature 
function $K$ for $S_g$ is then explicitly given by  
$$
 K(x) = - e^{-2u(x)}\Delta u(x)\, ,
\eqno\eqlbl\PDE
$$
where $\Delta$ is the Laplacian for the standard metric on $\RR^2$. 
The quantity 
$$
{\cal K}(u)\equiv \int_{\RR^2} K(x) e^{2u(x)}\dx \, ,
\eqno\eqlbl\totCdef
$$
where $\dx$ denotes Lebesgue measure on $\RR^2$, is called the 
integral curvature of the surface (sometimes called total curvature).
We say that $S_g$ is a {\it classical} surface over $\RR^2$
if $u \in C^2(\RR^2)$. Clearly, $K\in C^0(\RR^2)$ in that case. 

The inverse problem, namely to prescribe $K$ and to find a surface $S_g$
pointwise conformal to $\RR^2$ for which $K$ is the Gauss curvature,
renders (\PDE) a semi-linear elliptic PDE for the unknown function $u$. 
The problem of prescribing Gaussian curvature thus amounts to studying the 
existence, uniqueness or multiplicity, and classification of solutions $u$
of (\PDE) for the given $K$. A particularly interesting aspect of 
the classification problem is the question under which conditions  
radial symmetry of the  prescribed  Gauss curvature function $K$ 
implies radial symmetry of the classical surface $S_g = (\RR^2, g)$,
and under which conditions  radial symmetry is broken.

Notice that the inverse problem may not have a solution. In
particular, when considered on $\SS^2$ instead of $\RR^2$, there
are so many obstructions to finding a solution $u$ to 
(the analog of) (\PDE) for the prescribed $K$ that Nirenberg 
was prompted many years ago to raise the question:
``Which real-valued functions $K$ are Gauss curvatures of some 
surface $S_g$ over $\SS^2$?'' For Nirenberg's problem, see 
[\Beckner, \Branson, \CarlenLossB, \ChangYangA, \ChangYangB, \Han,
\KazdanWarner, \KiePHYa,  \LiYY,  \Ma, \MoserA, \MoserB, \Obata, \Onofri].
For related works on other compact 2-manifolds, see e.g.
[\DJLW,  \Tarantello].

In this work we are interested in the prescribed Gauss curvature
problem on $\RR^2$. There is a considerable literature on this 
problem, e.g.
[\Aubin, \Aviles, \ChaKieCMP, \ChengLinMA, \ChengLinJDE, \ChengLinASNSP,
\ChengNiI, \ChengNiII, \McOwen, \Ni, \Sattinger].
We here will study the existence problem of surfaces for a large class 
of $K$ via a novel approach. We also mention an existence-of-solutions result 
for a monotonically decreasing $K$ that is  unbounded below and positive at 
the origin. Moreover, we study the question of radial symmetry of classical 
surfaces, which correspond  to classical solutions of (\PDE), for 
monotonically decreasing $K$ and for non-positive $K$.

The problem of non-positive prescribed Gaussian curvature $K$
is already fairly well understood, see 
[\Ahlfors, \ChengNiI, \ChengNiII, \Osserman, \Poincare, \Wittich]. 
In particular, Theorem III of [\ChengNiI] 
characterizes any $S_g$ with compactly supported $K$ and finite integral 
curvature uniquely by its integral curvature and by an entire harmonic 
function $H$ to which $u$ is asymptotic at infinity. If the entire 
harmonic function is constant and $K$  radially symmetric, 
then $u$ is radially symmetric, by uniqueness. 
Theorem II of [\ChengNiI] characterizes any $S_g$ with 
$K\sim -C|x|^{-\ell}$ when $|x|\to\infty$,  $\ell>2$, and finite 
integral curvature uniquely by its integral curvature alone, so that
$u$ is radially symmetric if $K$ is. Theorem II of [\ChengNiI] is
extended in [\ChengNiII] to $K$ satisfying  an integrability
condition and $C|x|^{-m}\leq |K(x)|\leq C|x|^m$ as $|x|\to \infty$.

Our Theorem 2.1 below generalizes Theorem III of [\ChengNiI] 
as well as Cheng-Ni's Theorem II [\ChengNiI] and its sequel
in [\ChengNiII] to a larger class of $K$ satisfying mild integrability 
conditions without pointwise asymptotic bounds or even 
compact support for $K$. Our existence results follow as 
corollaries from our probabilistic Theorem 8.4 that applies to 
non-negative $K$ as well as non-positive ones. 
We prove our Theorem 8.4 in section VIII using the methods
developed in [\KieCPAM, \clmpCMPcan, \KieJllLMP] and [\KieSpo];
see also [\KiePHYa]. For non-positive radial $K$ the radial symmetry 
of $u$ then follows from our uniqueness Theorem~2.2, which we
prove in its dual version Theorem~9.1 in section IX.  

Prescribing Gaussian curvature $K$ which is somewhere strictly positive
is a much richer problem  and less well understood. Existence results
are available in 
[\Aviles, \ChengLinJDE, \ChengLinASNSP,  \McOwen, \Sattinger]; 
note [\ChengLinASNSP] regarding [\Aviles].
The question of radial symmetry of $u$ has been 
studied by various authors for decreasing $K$ 
under various additional conditions, see 
[\ChaKieCMP, \ChenLiA, \ChenLiB, \ChengLinMA, \PrajapatTarantello]. 

As already emphasized above, our Theorem 2.1 establishes existence of $u$ 
also for non-negative $K$, under mild integrability conditions on $K$ 
rather than prescribed asymptotic behavior or pointwise bounds as employed in 
[\Aviles, \ChengLinJDE, \ChengLinASNSP, \McOwen, \Sattinger]. 
We also announce an existence result of a radial surface with positive 
integral curvature for a radial continuous $K$ that is positive at 
the origin and diverges logarithmically to $-\infty$ as $|x|\to \infty$, 
see our Proposition 2.4. In our proof of Proposition 2.4 we actually do 
not prescribe $K$ but, inspired by [\KieJllPHPL], we consider a 
{\it system of equations} whose solutions determine both $K$ and $u$, 
and we use scattering theory  and gradient flow techniques to control it. 
This system case is of independent interest, and details of which will be 
published elsewhere.

The radial symmetry of surfaces with $K$ positive somewhere does
not follow simply by uniqueness. In section II we list various 
non-radial surfaces with radial Gauss curvature. We extract from 
this discussion a set of conditions on $K$ and $g$ which 
rule out the various non-radial surfaces we found. In particular, 
we demand $K$ be radial decreasing. We formulate a conjecture that under this 
set of conditions any classical surface for the corresponding prescribed 
Gauss curvature $K$ is radially symmetric about some point. 

We then state (section III), and subsequently prove (sections IV-VII),
using  the method of moving planes [\GidasNiNirbrg, \Ni], our 
Theorems 3.6 and 3.7 on radial symmetry of classical surfaces. 
Our symmetry theorems require a  slightly stronger set of conditions than 
formulated in our conjecture. However, our conditions are considerably weaker
than those used in the papers [\ChaKieCMP, \ChenLiA, \ChenLiB, \ChengLinMA].
In particular, we impose no pointwise bounds near infinity on positive 
$K$. We also allow $K$ to be unbounded below, but then with 
some growth conditions near infinity, allowing logarithmic as 
well as power law growth of $|K|$. Our existence-of-solutions 
Theorem 2.1 and Proposition 2.4 establish that solutions exist
under these conditions on $K$ and thus verify
that our radial symmetry theorems cover more cases than the
earlier symmetry results listed above. 

After submission of our work, existence results when $K$ is positive 
somewhere and satisfies $0\geq K(x)\geq -C|x|^{\ell}$ as $|x|\to\infty$, with 
$0<\ell <2$, appeared in [\ChengLinNA]. Of these surfaces, those which
also satisfy the hypotheses on $K$ listed in our Proposition 3.5 are 
radial symmetric by our Theorems 3.6 and 3.7.
 \vfill\eject

\medskip
\chno=2
\equno=0
\noindent
{\bf II. BROKEN SYMMETRY AND A SYMMETRY CONJECTURE}
\smallskip

We say that $S_g$ is radially symmetric about some point 
$x^\ast\in \RR^2$ if the associated solution $u$ of (\PDE) 
satisfies $u(x-x^\ast) = u(\cR(x-x^\ast))$ 
for any $\cR \in SO(2)$. We say that $u$ is non-radial if no such 
point exists. We now collect a list of examples of non-radial
surfaces from which we extract conditions on $K$ and $g$
under which one can hope to assert the radial symmetry of $u$.

Clearly, $u$ cannot be radially symmetric about some point if
$K$ is not radially symmetric about the same point. Without loss 
we choose the point about which $K$ is radially symmetric 
to be the origin, i.e. we demand
$$
K(x) = K(\cR x)\, .
\eqno\eqlbl\Krotinv
$$
A few moments of reflection reveal that some further conditions 
on $K(x)$ and $u(x)$ will be needed, for without further conditions, 
examples to non-radially symmetric surfaces having a Gauss curvature 
$K$ satisfying (\Krotinv) are readily found.

In particular, if $K$ satisfying (\Krotinv) is compactly supported 
then solutions $u$ of (\PDE) that display some non-constant entire harmonic 
behavior near infinity have been asserted to exist (for non-positive $K$) 
in Theorem III of [\ChengNiI]. Our first theorem, proved in section VIII,
generalizes Theorem III of [\ChengNiI], as well as their Theorem II
and its extension in [\ChengNiII], to a much wider class of 
sufficiently `concentrated' $K$ that have well defined sign.
We define the sign  $\sigma (K)$ of the function $K$ by:
$\sigma(K) = +1$ if $K\not\equiv 0,\ K(x)\geq 0$ for all $x\in \RR^2$; 
$\sigma(K) = -1$ if $K\not\equiv 0,\ K(x)\leq 0$ for all $x\in \RR^2$;  
$\sigma(K) = 0$ if $K(x)\equiv 0$. For  other $K$, $\sigma(K)$ does not exist.

\smallskip
\noindent
{\bf Theorem 2.1:} {\it 
		Assume $K\in L^\infty(\RR^2)$ has well defined
		sign $\sigma (K)$. 
		Assume furthermore that for some entire harmonic function 
		$H:\RR^2\mapsto \RR$, and all $0 < \gamma < 2$, $K$ satisfies 
$$
\int_{{\rm B}_1(y)}
|y - x|^{-\gamma} |K(x)|e^{2H(x)}\dx \longrightarrow 0
\qquad	{\rm as}\ |y| \to \infty\, , 
\eqno\eqlbl\KHgamma
$$
where ${\rm B}_R(y)\subset \RR^2$ is the open ball of radius $R$ 
centered at $y$. Given the same $H$, assume also that $K$ satisfies
$$
\int_{\RR^2} |K(x)| e^{2H(x)}|x|^{q} \dx < \infty\, 
\eqno\eqlbl\KHq
$$
		for some $q >0$. If $K\leq 0$, define
$$
\kappa^*(K,H) = - 2\pi\,
\sup_{q>0}\, \bigl\{ q :\ {\rm (\KHq)\ is\ true\, } \bigr\} \, .
 \eqno\eqlbl\kappaSUPstar
$$   
		Then, for any such $K$, $H$, and any $\kappa$ satisfying 
$$
\kappa \in \cases{
\ (\kappa^*,0)		&{\rm if} $K\not\equiv 0,\ K\leq 0$;	\cr
\qquad     \{0\}	&{\rm if} $K\equiv 0$;			\cr
\qquad     (0,4\pi)  	&{\rm if} $K\not\equiv 0,\ K\geq 0$,	\cr }
\eqno\eqlbl\kappacases
$$
		there exists a solution 
		$u = U_{H,\kappa}\in W^{2,p}_{\rm loc}\cap L^\infty_{\rm loc}$ 
		of (\PDE) for the prescribed Gaussian curvature function $K$, 
		having integral curvature
$$
{\cal K}\bigl(U_{H,\kappa}\bigr) ={\kappa}\, ,
\eqno\eqlbl\totCH
$$
		and having asymptotic behavior given by
$$
\qquad
U_{H,\kappa}(x) = H(x) - {\kappa \over 2\pi} \ln |x| + o(|\ln |x||)  
\qquad {\rm as}\ |x| \to \infty\, .
\eqno\eqlbl\Uasymp
$$

	If moreover $K\in C^{0,\alpha}$, then $U_{H,\kappa}$ is
a classical solution.
	If that $K\in C^{0,\alpha}$ also satisfies (\Krotinv), and $H$ 
is non-constant, 
then $U_{H,\kappa}$ generates a classical surface which is asymptotic to 
a non-radial entire harmonic surface, hence breaking radial symmetry.
}\smallskip

We remark that, if $|K|\in C^{0,\alpha}$ satisfying (\Krotinv) is 
also decreasing, then all these conclusions hold without imposing (\KHgamma).

Surfaces which are asymptotic to some non-radial entire harmonic 
surface (entire harmonic surfaces for $K\equiv 0$) 
can be eliminated by the mild integrability condition 
$$
u^+ \in L^1\left({\rm B}_R(y),\dx\right),\ \ {\rm uniformly\ in}\ y\, ,
\eqno\eqlbl\uplusLone
$$
where $u^+(x) := \max\{u(x),0\}$. 
For the category of $K\leq 0$ covered in Theorem 2.1 which
in addition satisfy (\Krotinv), condition (\uplusLone) already 
eliminates {\it all} non-radial solutions $u$ of (\PDE) with finite 
integral curvature. Indeed, we have,

\smallskip
\noindent
{\bf Theorem 2.2:} {\it Under the hypotheses stated in Theorem 2.1 
and in (\uplusLone), if $K \leq 0$, then the solution $U_{H,\kappa}$ is 
unique. Moreover, if $K\leq 0$ also satisfies (\Krotinv), then 
$U_{H,\kappa}$ is radial symmetric and decreasing.}
\smallskip

It remains to discuss $K$ which are strictly positive somewhere.
In that case, among the $S_g$ that satisfy (\Krotinv) and (\uplusLone)
one finds non-radial surfaces that are periodic about the origin 
of the Euclidean plane, having fundamental period $2\pi/n$, with $n>1$. 
We illustrate this with the following examples,
taken from [\ChaKieCMP] (see also [\PrajapatTarantello]). For $x\neq 0$, 
we introduce the usual polar coordinates $(r,\theta)$ of $x$, i.e. 
$r=|x| >0$ and $\tan\theta = x_2/x_1$, with $\theta\in [0,2\pi)$. 
Let $\NN$ denote the natural numbers. For $n\in\NN$, let 
$K(x)=K^{(n)}(x)$, with
$$
K^{(n)}(x) = 4n^2 |x|^{2(n-1)}\, .
\eqno\eqlbl\chakieK
$$ 
Clearly, $K^{(n)}\in C^\infty(\RR^2)$. Let $y \in\RR^2$
be chosen arbitrarily, except that $y \neq 0$, and let $\theta_0$
be the polar angle coordinate of $y$. Let $\zeta\in \RR$. 
Then $u(\, .\,) = U_{\zeta}^{(n)}(\, .\,;y)$, with 
$$
\eqalignno{
U_{\zeta}^{(n)}(x;y) =
&  
  - \ln \left( 1 - 2{|x|^n\over |y|^n}\cos\bigl(n(\theta-\theta_0)\bigr)
	\tanh\zeta + {|x|^{2n}\over |y|^{2n}}\right)
&\cr &
 - \ln \Bigl(|y|^{n}\cosh \zeta \Bigr)
&\eqlbl\chakieSOL\cr}
$$
is a $C^\infty(\RR^2)$ solution of (\PDE) for the Gaussian
curvature function (\chakieK). The integral curvature of  the
surface described by (\chakieSOL) is given by   
$$
{\cal K}\bigl(U_{\zeta}^{(n)}(x;y)\bigr) = 4\pi n\, ,
\eqno\eqlbl\chakietC
$$
independently of $\zeta$ and $y$. For $\zeta = 0$ and all $n\in\NN$,
the solution (\chakieSOL) is radially symmetric about the origin. 
For $\zeta \neq 0$, if $n =1$ so that (\chakieK) reduces 
to a constant, $K^{(1)}= 4$, the solution (\chakieSOL) is periodic 
about the origin with fundamental period $2\pi$, yet it is radially 
symmetric about and decreasing away from the point 
$x^\ast =\tanh(\zeta)y$. For $\zeta\neq 0$ and $n>1$,
in which cases $K^{(n)}$ increases monotonically with $|x|$,  
the solution (\chakieSOL) is periodic about the origin 
with fundamental period $2\pi/n$, whence non-radial about any point;
see Figure~1.  

\medskip
\epsfxsize=8.5cm
\centerline{\epsffile{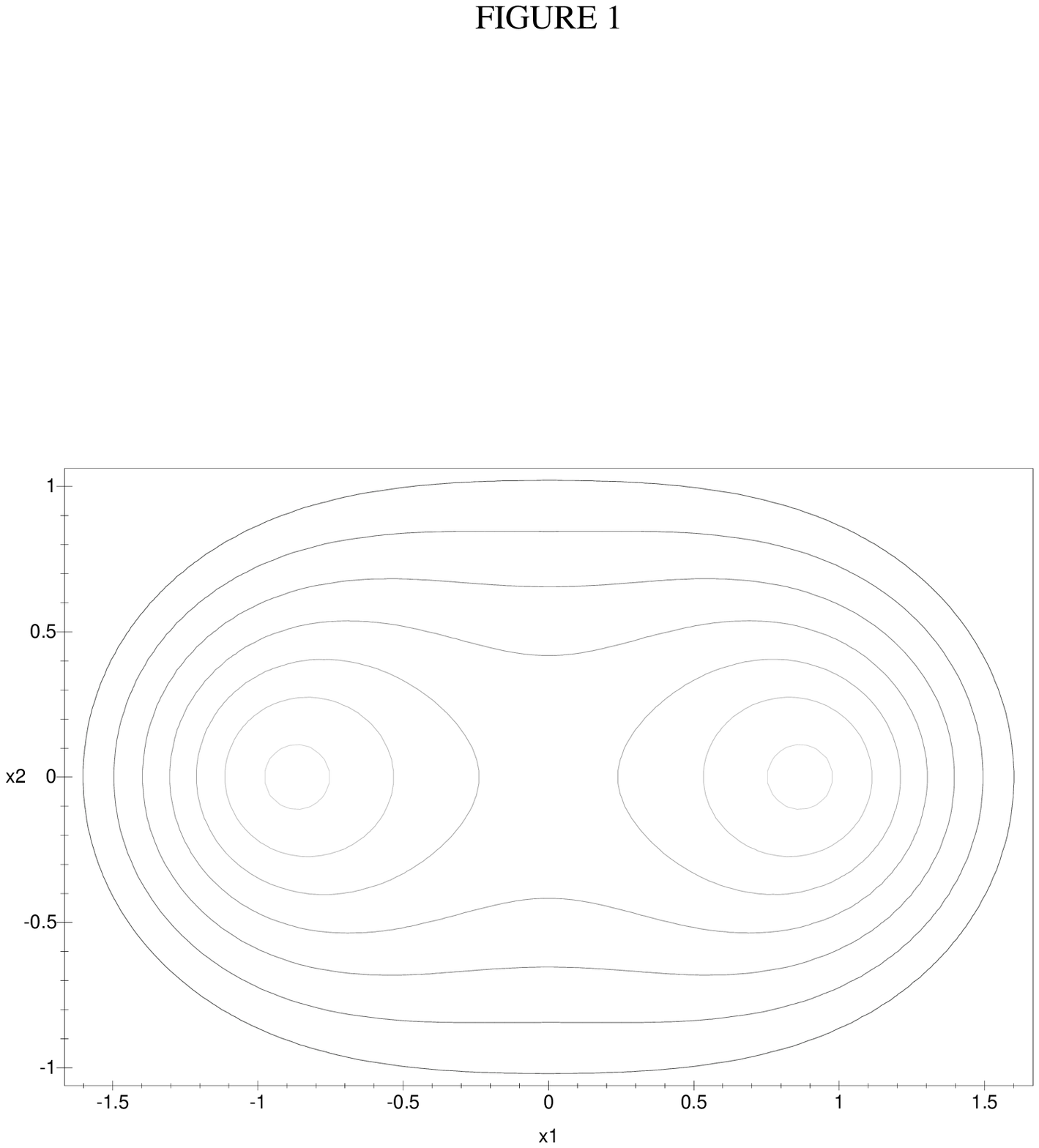}}

\centerline{
Fig.1:\ Level curves  $e^{2u(x)}= 2^a$, $a\in\{-5,-4,...,0,1\}$,
for $u$ given by (\chakieSOL)}

\centerline{
with $n=2$, $|y|=1$, $\theta_0=0$, $\zeta =1$. 
$\max e^{2u}\approx 2.57$ is taken at the centers of}

\centerline{
 the two islands. For $|x|$ large, the conformal factor 
 $e^{2u(x)} \sim C |x|^{-8}$  and the}

\centerline{
level curves become circular.$\phantom{more 
and more *******************}$}

\medskip

This last family of non-radial surfaces is eliminated by admitting
only monotonically decreasing radial $K$, i.e., those $K$ satisfying
$$
K(x) \leq K(y) \ \ \ {\rm whenever}\ \ \ |x| \geq  |y|\, .
\eqno\eqlbl\Kraddec
$$

Among the $S_g$ that satisfy (\Krotinv), (\uplusLone), and (\Kraddec),
we still find non-radial surfaces, namely when $K(x) = K_0$, with
$$
K_0 = \ {\rm constant}\, >0\, ,
\eqno\eqlbl\constK
$$
in which case (\PDE) is the conformally invariant Liouville equation 
[\Liouville]. Beside the radial symmetric entire solutions obtained 
with $n=1$ in (\chakieSOL), this equation has entire classical 
solutions that are periodic along a Cartesian coordinate direction. 
Let $y\in\RR^2$ be an arbitrary fixed point, and let 
$v\in \RR^2$ and $v^\prime \in \RR^2$ be two fixed 
vectors that are orthogonal w.r.t. Euclidean inner product, i.e. 
$\langle v,v^\prime\rangle =0$, having identical lengths given by 
$|v| = |v^\pr| = K_0^{1/2}$. Let $\zeta\in\RR$. 
Then  $u(\, .\, ) = U_\zeta(\, .\,;y)$, with 
$$
U_\zeta (x;y)
 = - \ln \bigl( \cosh(\zeta)\cosh \langle v,     x-y\rangle
	     - \sinh(\zeta)\sin  \langle v^\pr, x-y\rangle \bigr) \, ,
\eqno\eqlbl\stuartSOL
$$
is a non-radial $C^\infty(\RR^2)$ solution of (\PDE) for the 
Gauss curvature function (\constK); see also [\ChaKieCMP]. 
For $\zeta =0$, the 
solution is translation invariant along $v^\pr$, while for 
$\zeta\neq 0$ it is periodic along $v^\pr$ with period
$2\pi/\sqrt{K_0}\,$, see Figure~2. 

\smallskip
\epsfxsize=8.5cm
\centerline{\epsffile{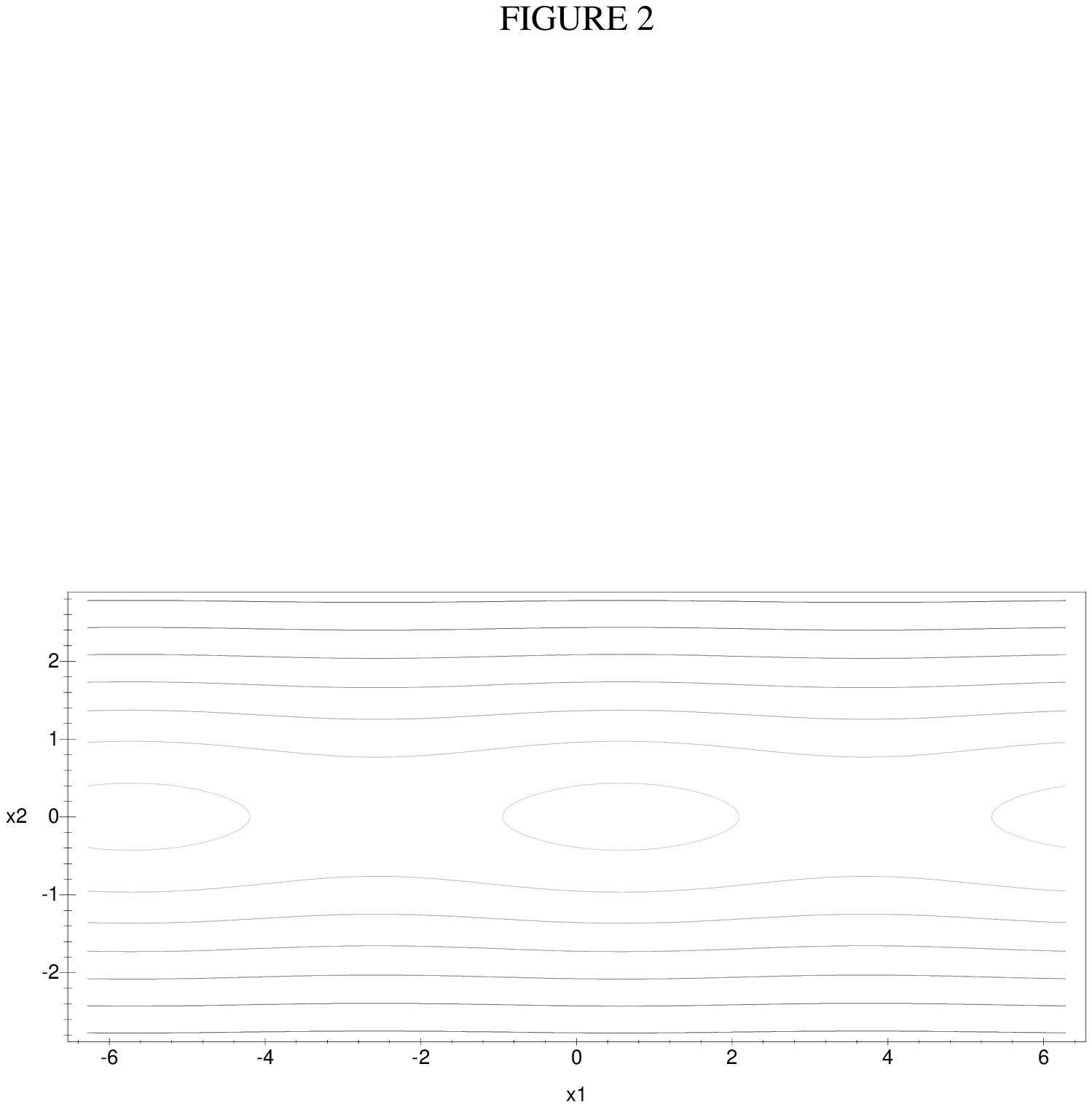}}

\centerline{
Fig.2: Level curves $e^{2u}= 2^{a}$, $a\in\{-6,-5,...,0\}$, with $u$ 
given by}
\centerline{
(\stuartSOL), with $\zeta =1$, $y=- v^\prime$, $K_0=1$, 
$x_1=\langle x,v^\pr\rangle$ and $x_2=\langle x,v\rangle$.}

\centerline{ 
$\max e^{2u}\approx 1.22$ is taken at the centers of the islands.
For $|\langle v,x\rangle|$ }

\centerline{
large, $e^{2u(x)} \sim C e^{-|\langle v,x\rangle|}$ and level 
curves become straight lines.}
\vfill\eject

Since $\exp(2U_\zeta(\, .\,;y))\not\in L^p(\RR^2, \dx)$ for all 
$p$ except $p=\infty$, the surface corresponding to (\stuartSOL) 
has integral curvature ${\cal K}(u) =+\infty$, as does any surface 
that is periodic or invariant along a fixed direction. 

To rule out translation invariant surfaces and 
those that are periodic along a fixed Cartesian direction of the 
Euclidean plane, we could impose the integrability condition
$\int \exp(2u(x)) \dx < \infty$. However, it suffices to 
impose the milder, and more natural, restriction that 
the surface's Gauss curvature is absolutely integrable, i.e.
$$
\int_{\RR^2} |K(x)| e^{2u(x)} \dx < \infty\, ,
\eqno\eqlbl\tCbound
$$
which reduces to $\int \exp(2u(x)) \dx < \infty$ if $K= const. >0$.
%\vfill\eject

We summarize the various conditions on $S_g$ as follows.

\smallskip
\noindent
{\bf Definition 2.3:} {\it 
		For each $K \in C^{0,\alpha}(\RR^2)$ satisfying (\Kraddec), 
		we denote by ${\bf S}_K$ the set of classical surfaces 
		$S_g$ with Gauss curvature $K$ being absolutely integrable,
		(\tCbound), and with metric (\metric) satisfying (\uplusLone).
}
\smallskip

Notice that there exist $K$ for which the set ${\bf S}_K$ is empty. 
Thus, since $K$ satisfies (\Kraddec), no  entire solutions of (\PDE)
exist if $K <0$ everywhere [\Osserman]. In particular, 
entire solutions in $\RR^2$ with $K= constant\ <0$ do not exist, 
see [\Ahlfors, \Osserman, \Wittich]. Moreover, 
if $K(x)\sim - C|x|^p$ for $p\geq 2$ (irrespective of whether
$K(x) \leq 0$ for $|x|<R$ or not) then it follows from an easy
application of Pokhozaev's identity that ${\bf S}_K$ is empty.

On the other hand, if $K\geq 0$ everywhere, then there are plenty of 
radially symmetric surfaces in ${\bf S}_K$, which follows
from our Theorem 2.1 with $H\equiv constant$.
Furthermore, we note that ${\bf S}_K$ is not
empty for certain radial $K$ that are unbounded 
below, for

\smallskip
\noindent
{\bf Proposition 2.4:} {\it There exist continuous $K(x)$ satisfying 
	(\Kraddec) and $K(x)\sim -C\ln |x|$ as $|x|\to \infty$
	for which ${\bf S}_K$ contains radial surfaces with 
	finite positive integral curvature.  
}
\smallskip
The proof of Proposition 2.4, which uses ideas from scattering theory
similar to those in [\KieJllPHPL] together with gradient flow techniques,
is of independent interest and will be published elsewhere. 

All known examples of surfaces in ${\bf S}_K$ are radially symmetric,
and we could not conceive of any counterexample to radial symmetry.
Hence, we conjecture that all surfaces in ${\bf S}_K$ are radially 
symmetric. More precisely, our conjecture reads as follows.
\smallskip
\noindent
{\bf Conjecture 2.5:} {\it 
			Any classical surface $S_g \in {\bf S}_K$ 
			is equipped with a radially symmetric non-expansive 
			metric, in the sense that the conformal	factor 
			$e^{2u}$ is radially symmetric and decreasing 
			about some point.}
\medskip

Presumably, Conjecture 2.5 can even be widened to include certain $K$
that are not everywhere decreasing, see [\ChaKieCMP, \PrajapatTarantello] 
for examples.
However, currently it seems not  clear how to prove even Conjecture 2.5 
without some additional technical conditions. In the ensuing sections 
we will first state and then prove radial symmetry theorems for 
${\bf S}_K$ under conditions that are weaker than those used in 
previous theorems, yet slightly stronger than those stated
in Conjecture 2.5. In the next section we state precisely
our main symmetry results, assess the territory covered by them, 
and also compare them to existing results. 

 \vfill\eject

\medskip
\chno=3
\equno=0
\noindent
{\bf III. SYMMETRY THEOREMS FOR RADIAL DECREASING $K$}
\smallskip

To state our new symmetry results for $K \in C^{0,\alpha}(\RR^2)$ satisfying 
(\Kraddec), we define
$$
\kappa_*(K) = \pi\inf \left\{ q >0  :\, 
\int_{\RR^2} |K(x)| (1+|x|)^{-q}\dx < \infty \right\}\, .
\eqno\eqlbl\kappastar
$$
The significance of $\kappa_*(K)$ is that of 
an explicit lower bound to the integral curvature. 

\smallskip
\noindent
{\bf Proposition 3.5:} {\it Let $K \in C^{0,\alpha}(\RR^2)$ satisfy (\Kraddec).
	If $K$ is unbounded below, then let $K$ also satisfy one of the 
	following two conditions, either (1): there exists some $C >0$ 
	such that
$$
 |K(x)| \leq C \inf_{y\in {\rm B}_1(x)}  |K(y)| 
\ \ {\rm as}\ \ |x|\to
\infty\, ,
\eqno\eqlbl\kappastarCone
$$
	uniformly in $x$ (this condition is satisfied, e.g.,
	if $K\sim -C|x|^\ell$, any $\ell >0$); 
	or (2): there exist some finite $P\geq 1$ and $C>0$ 
	such that
$$
|K(x)| \leq C |\ln |x||^P\ \ {\rm as}\ |x|\to\infty\, .
\eqno\eqlbl\kappastarCtwo
$$
Let $K$ be the Gauss curvature function for
	a surface $S_g \in {\bf S}_K$. Then the integral curvature 
	of $S_g$ is bounded below by
$$
{\cal K}(u) \geq \kappa_*(K) \, .
\eqno\eqlbl\tCkappastar
$$}

We now state two theorems on radial symmetry of surfaces in 
${\bf S}_K$, distinguishing the cases ${\cal K}(u) > \kappa_*(K)$ 
and ${\cal K}(u) = \kappa_*(K)$. Our first theorem (Theorem 3.6 below)
verifies Conjecture 2.5, under the hypotheses of Proposition 3.5, 
for all integral curvatures ${\cal K}(u) > \kappa_*(K)$. 
By Proposition 3.5, this covers the spectrum of potential integral 
curvature values all the way down to its lower bound (\tCkappastar), 
but not including it. This signals that the borderline case 
${\cal K}(u) = \kappa_*(K)$ is  critical. 
The critical case ${\cal K}(u) = \kappa_*(K)$ is dealt with
in our Theorem 3.7 below, were we assert the radial symmetry 
and decrease of $u$ under an additional hypothesis which is 
mildly stronger than (\tCbound). 

\smallskip
\noindent
{\bf Theorem 3.6  (the sub-critical case):} {\it Under the
		assumptions stated in Proposition 3.5, all 
		surfaces $S_g \in {\bf S}_K$ with integral 
		curvature ${\cal K}(u) > \kappa_*(K)$ are equipped with a 
		radially symmetric, non-expansive metric (\metric), i.e. 
		there exists a point $x^\ast\in \RR^2$ such that $u$ in 
		(\metric) is radially symmetric and decreasing about $x^\ast$, 
$$
u(x-x^\ast)\leq u(y-x^\ast)\ \ \ {\rm whenever}\ |x-x^\ast| \geq |y-x^\ast|\, .
\eqno\eqlbl\uraddec
$$
		Moreover, if $K \not\equiv \, constant$, then $x^\ast = 0$, 
		and if $K\equiv \, constant$, then $x^\ast$ is arbitrary.} 
\smallskip
\noindent
{\bf Theorem 3.7  (the critical case):} {\it 
		Under the assumptions stated in Proposition 3.5, a surface 
		$S_g \in {\bf S}_K$ having integral curvature
		${\cal K}(u) = \kappa_*(K)$ is equipped with a 
		radially symmetric, non-expansive metric (\metric) (in 
		the sense of (\uraddec)) provided
$$
\int_{\RR^2} \bigl|\ln |x|\bigr|^2 |K(x)|\, e^{2u(x)} \dx < \infty\, .
\eqno\eqlbl\intKeulnsq
$$
		In that case, if $K \not\equiv \, 0$, then $x^\ast = 0$, 
		and if $K\equiv \, 0$, then $x^\ast$ is arbitrary.
}
 \vfill\eject

With reference to Conjecture 2.5, the foremost question now is 
how much of ${\bf S}_K$ is actually covered by our Theorems 3.6 
and 3.7, and how much remains uncharted territory. A priori speaking, 
Theorems 3.6 and 3.7 leave us anywhere in between the following extreme 
scenarios. In the best conceivable case, all surfaces with critical 
integral curvature satisfy (\intKeulnsq), and then Theorems 3.6 and 3.7 
taken together would prove Conjecture 2.5 completely. In the worst 
conceivable case, all surfaces have critical integral
curvature, and none satisfies (\intKeulnsq), in which case
Theorems 3.6 and 3.7 would be empty. 
		% together would not prove Conjecture 2.5. 
To assess the situation, we 
need to address the question whether for any $K$ 
there exists a critical surface $S_g$ such that
inequality (\tCkappastar) is an equality, and if, 
whether any such critical $S_g$ satisfies (\intKeulnsq).
Notice that (\intKeulnsq) is only needed for those $K$ for 
which there exists a critical surface, i.e. a surface for which
(\tCkappastar) is an equality. 

Inequality (\tCkappastar) is certainly an equality in the trivial 
case $K\equiv 0$, where we have ${\cal K}(u) = 0 = \kappa_*(0)$.
Of course, (\intKeulnsq) is trivially satisfied when  $K\equiv 0$, 
whence this case is covered by Theorem 3.7. 

If $K(x)\not\equiv 0$ decreases to zero at least as
$C|x|^{-2-\eps}$, possibly having compact support, 
then ${\cal K}(u) > 0$, by (\totCdef),  while $\kappa_*(K)= 0$.
Obviously inequality (\tCkappastar) is strict in these cases, 
whence Theorem 3.6 covers all possible surfaces for each such $K$. 
We remark that by a our Theorem 2.1 with 
$H\equiv constant$ it follows for such 
decreasing $K$ that surfaces do exist for all integral 
curvature values in the open interval $(0,4\pi)$. Together 
with ${\cal K}(u)>0$, this implies for these $K$ that 
$\kappa_*(K)= 0$ is the infimum to the set of integral curvatures 
for surfaces $S_g \in {\bf S}_K$.

For Gauss curvature functions $K=K_0 \, >0$, with 
$K_0$ a constant, we have $\kappa_*(K_0) = 2\pi$, 
while ${\cal K}(u) = 4\pi$ for all solutions $u$ of (\PDE), (\uplusLone), 
(\Kraddec), (\tCbound), see [\ChaKieGAFA, \ChenLiA]. Not only  is
inequality (\tCkappastar) strict in these cases,  
$\kappa_*(K_0)$ is not even the best constant in the sense of an
optimal lower bound to the integral curvature. Clearly, the cases 
$K=K_0 \, >0$, with $K_0$ a constant, are entirely covered by Theorem 3.6. 

The situation seems less clear when, as $|x|\to \infty$,  
$K$ behaves like $C|x|^{-p}$ or like $-C |x|^p$, with $p<2$. 
In these cases, explicit existence statements of surfaces 
in ${\bf S}_K$ with critical curvature ${\cal K}(u) = \kappa_*(K)$ 
seem currently not available. 

We remark that surfaces with critical curvature ${\cal K}(u) =
\kappa_*(K)$ do exist 
when $K\leq 0$ and $K(x)\sim -|x|^{-\ell}$ as $|x|\to \infty$, with 
$\ell >2$. While those surfaces are radially symmetric by a uniqueness 
argument, it is nevertheless quite interesting to register that 
they do {\it not} satisfy (\intKeulnsq)! The metric (\metric) of 
these surfaces is equipped with a conformal factor $e^{2U}$, where 
$U$ is the maximal solution of Cheng and Ni [\ChengNiI], see their 
Theorem II, p. 723. Cheng and Ni's result signals the possible 
existence of surfaces with critical curvature in ${\bf S}_K$ to
which our Theorem 3.7 does not apply.

We summarize this state of affairs with the following list of
interesting open questions. 
\smallskip
\noindent
{\bf Open Problems 3.8:} {\it Do there exist
radially decreasing $K\not\equiv 0$ for which there exist
solutions of (\PDE), (\uplusLone), (\Kraddec) with 
${\cal K}(u) = \kappa_*(K)$? If the answer to the previous 
question is positive, is (\intKeulnsq) a genuine condition, in the 
sense that there exist surfaces in ${\bf S}_K$
violating (\intKeulnsq)? And in case the answer to that question
is also positive, is Conjecture 2.5 false for some of these surfaces? }
\smallskip

Incidentally, the above discussion also points to a 
related open question which, though less directly relevant to our 
inquiry into radial symmetry, is an interesting problem in itself. 
To this extent, we introduce the notion of a least integrally curved 
surface in ${\bf  S}_K$, and, with an eye toward the above discussion,
also the notion of when such a surface is critical.

\smallskip
\noindent
{\bf Definition 3.9:} {\it 
		A surface $S_g\in {\bf S}_K$ is called 
		{\bf least integrally curved} if 
		${\cal K}(u)= \kappa(K)$, where $\kappa(K)$ is 
		defined as the infimum of the set of integral 
		curvatures for which there exists a surface 
		$S_g\in {\bf S}_K$, given $K$. A least
		integrally curved surface is called {\bf critical}
		if $\kappa(K) = \kappa_*(K)$.
}

\smallskip
\noindent
{\bf Open Problems 3.10:} {\it Find and classify all $K$
for which there exists a least integrally curved surface 
in ${\bf S}_K$! With reference to Problems 3.8, 
determine which of those surfaces are critical!
} 
\smallskip

We now return to the question of radial symmetry and to our strategy 
of proof for our Theorems 3.6 and 3.7. We use the technique of the 
moving planes [\GidasNiNirbrg, \LiCM], adapted to the setting in 
two-dimensional Euclidean space (where it is proper to rather speak 
of moving lines) so that it is possible to move in the lines from 
`spatial infinity.'
Due to the logarithmic divergence of solutions at infinity, 
this part is more delicate than in higher dimensions, in 
particular when $K>0$. Various authors before have applied this 
method  to the problem under consideration here. Hence, 
before we enter the details of our proof, we briefly explain 
in which way our Theorems 3.6 and 3.7 go beyond existing results. 
 
Radial symmetry of surfaces with strictly positive, {\it constant}
Gauss curvature function (\constK) and finite integral curvature 
(\tCbound) was proven by Chen and Li [\ChenLiA]. 
In [\ChenLiA], a radial `comparison function' was invented that 
made it possible to overcome the `problem at infinity.' 
In this case the result allows one to compute all surfaces explicitly,
which are given by (\chakieSOL) with $n=1$. This result was also 
obtained, with two different alternate methods, in [\ChouWan] and 
in [\ChaKieGAFA]. 

In [\ChaKieCMP] the method of [\ChenLiA] was extended to a  wider 
class of surfaces with monotone decreasing, bounded Gauss 
curvature functions, given certain integrability conditions.  
The following was proven in [\ChaKieCMP].

% \vfill\eject

\smallskip
\noindent
{\bf Theorem 3.11:} {\it 
		Let $K$ be the bounded Gauss curvature function of a 
		classical surface $S_g$, with metric given by 
		(\metric), and assume that  (\uplusLone), (\tCbound), 
		and (\Kraddec) are satisfied. Let $K^+$ denote the 
		positive part of $K$. Then any  surface $S_g$ whose
		integral curvature satisfies
$$
{\cal K}(u) > 
\pi\biggl( 3 +  \limsup_{|x|\to\infty}\, {\ln K^+(x)\over \ln |x|}\biggr) \, ,
\eqno\eqlbl\threePIlimsup
$$
		is radial, more precisely there exists a point 
		$x^\ast\in \RR^2$ such that (\uraddec) holds. 
}

\smallskip
\noindent
{\bf Remark 3.12:} {\it 
		The proof of Theorem 3.11 is contained 
		in [\ChaKieCMP], proof of Theorem P1.}
\smallskip

Clearly, Theorem 3.11 falls short of proving Conjecture 2.5, for one 
because $K$ is assumed bounded in Theorem 3.11, and furthermore because
there exist surfaces with radial decreasing and bounded 
Gauss curvature function whose integral curvatures ${\cal K}(u)$ 
violate (\threePIlimsup). For instance, consider the special case 
of (\KHq) where $K>0$ satisfies the growth condition
$$
\lim_{|x|\to\infty}\, {\ln K(x)\over \ln |x|} = -m < -2\, .
\eqno\eqlbl\Kgrowth
$$
Theorem 3.11 asserts the radial symmetry of surfaces
with ${\cal K}(u) >\pi(3-m)^+$ ((\threePIlimsup) 
with `$\limsup$' now `$\lim$'). Surfaces with integral curvature 
in the interval $0 < {\cal K}(u) \leq \pi (3-m)^+$, which 
by our Theorem 2.1 exist for $m \in (2,3)$, are not covered 
by Theorem~3.11.

On the other hand, by Proposition 3.5, $\kappa_*(K)=0$ for
$K>0$ satisfying (\Kraddec) and (\KHq), while 
${\cal K}(u) > 0$ because of $K>0$. Hence, our Theorem 3.6 
applies and asserts the radial symmetry of all surfaces in 
${\bf S}_K$ with non-negative radially decreasing Gauss 
curvature functions $K$ satisfying (\KHq), including  
as special case the $K$ that satisfy (\Kgrowth). 

Closer inspection of the proof of Theorem 3.11 (see Remark 3.12) 
reveals that the origin of the $3\pi$ in (\threePIlimsup) versus the
$2\pi$ that is required to cover all surfaces for the $K$ satisfying 
(\Kgrowth) traces back to our using the comparison function 
of [\ChenLiA]. That comparison function, while well suited for 
constant and for certain monotonically decreasing Gauss curvature 
functions $K$, does not suit radial decreasing $K$ in general. 

One main technical innovation of the present paper is  the 
systematic construction of a new, radial comparison function which 
proves itself nearly optimal for handling the problem at infinity. 
We also obtain better control of solutions $u$ of (\PDE) near infinity, 
which allows us to forgo some technical contraptions used in [\ChaKieCMP]. 

Other, heuristic, comparison functions have been explored 
in the literature. Chen and Li in [\ChenLiB] use
a translation invariant comparison function rather than a radial one, 
and require the stronger conditions that $e^{2u}\in L^1(\RR^2)$,
thereby restricting integral curvatures to ${\cal K}(u) > 2\pi$,
to prove that all corresponding surfaces with
strictly positive, radially symmetric decreasing $K$ 
are given by radially symmetric and decreasing solutions $u$ of (\PDE). 
This result of [\ChenLiB] is contained in our Theorems 3.6 and 3.7.
Furthermore, it intersects with, but does not subsume, due to its stronger 
conditions on $u$,  Theorem 3.11 of [\ChaKieCMP]. For example, 
consider the Gauss curvature function $K(x) = K_\gamma(x)$, with
$$
K_\gamma(x) = 4 \gamma \exp\left(2(1-\gamma )U_0^{(1)}(x;y)\right) \, ,
\eqno\eqlbl\specialK
$$
where $U_0^{(1)}(x;y)$ is the special case $\zeta =0$ and $n=1$ in 
(\chakieSOL), with $y\neq 0$ arbitrary, and $0< \gamma \leq 1$. All 
$K_\gamma$ are radially decreasing, and we have 
$$
K_\gamma(x) \sim C |x|^{-4(1-\gamma)}\, .
$$
Clearly, 
$$
u(x) =  \gamma U_0^{(1)}(x;y)
\eqno\eqlbl\specialu
$$
is a radial, decreasing solution of (\PDE) for $K$ given  by (\specialK).
A classical radial surface described by (\specialu) has integral curvature
$$
\int_{\RR^2} K_\gamma(x) e^{2 \gamma U_0^{(1)}(x;y)}\dx
 = \gamma\, 4\pi \in (0,4\pi]\, ,
\eqno\eqlbl\specialtotC
$$
independently of $y$. When $\gamma \leq 1/2$, our examples (\specialu) 
violate Chen-Li's condition that $e^{2u}\in L^1$. 
Nevertheless,  for $K$ given by (\specialK), solutions of (\PDE) that 
satisfy (\uplusLone) and (\tCbound) also satisfy condition 
(\threePIlimsup) in Theorem 3.11, irrespective of $\gamma$, whence 
radial symmetry follows by Theorem 3.11. (Cf., also Theorem V2 in 
[\ChaKieCMP].) Incidentally, $\kappa_*(K_\gamma) =  2\pi(2\gamma -1)^+
< \gamma 4\pi$, and so none of these surfaces is critical. Hence, 
the radial symmetry of these surfaces follows by our Theorem 3.6 as well. 
Finally, a non-symmetric comparison function 
(a sum of a radial and a translation invariant function) is  
used in [\ChengLinMA] to prove the radial symmetry of surfaces with 
radial decreasing Gauss curvature function $K$, having 
finite integral curvature, under stronger conditions on $K$
than in our Theorems 3.6 and 3.7, namely that $K$ be strictly 
positive and decay slower than exponentially. 

This concludes our discussion of the radial symmetry theorems.
The next three sections of our paper are devoted to the proof of 
Theorems 3.6 and 3.7. In section VIII we prove Theorem~2.1, and 
in section IX we prove Theorem~2.2.

% \vfill\eject
 
\medskip
\chno=4
\equno=0
\noindent
{\bf IV. ASYMPTOTICS} 
\smallskip

To prepare the proofs of our theorems we need to gather some facts about
the asymptotic behavior of the solutions $u$ of (\PDE). In the following, 
$X = x/|x|^2$ denotes the Kelvin transform of $x$.

\smallskip
\noindent
{\bf Lemma 4.1:} {\it 
		Let $u$ be a classical solution of (\PDE) satisfying 
		(\uplusLone). Assume (\tCbound) holds, and that $K$ 
		satisfies (\Kraddec). Then $u$ satisfies the integral equation
$$
u(x) - u(0) = - {{\cal K}(u)\over 2\pi} \ln|x| - {1\over 2\pi} 
\int_{\RR^2} \ln | X - Y| K(y)e^{2u(y)} \dy
\eqno\eqlbl\uIntEQ
$$
		for all $x$.}

\medskip
{\it Proof:} By hypothesis, $K$ is monotone decreasing. We 
distinguish the cases with $K\geq 0$ from those where $K<0$ for $|x|>R$. 

In case $K$ becomes negative somewhere, say for $|x|>R$,  then outside
the disk ${\rm B}_R(0)$  the function $u$ is sub-harmonic, and so is $u^+$.  
Hence, for concentric  disks ${\rm B}_{1/2}(y)$ and ${\rm B}_1(y)$ 
we have 
$$
\left\|u^+\right\|_{L^\infty \left({\rm B}_{1/2}(y)\right)}\leq 
C \left\|u^+\right\|_{L^1\left({\rm B}_{1}(y)\right)}
\eqno\eqlbl\BhalfBone
$$
for some constant $C$ which is independent of $y$. Our hypothesis (\uplusLone) 
guarantees that the right side in (\BhalfBone) is bounded by a constant, 
whence we have a uniform $L^\infty$ bound for $u^+$ outside a disk, 
and this implies a uniform $L^\infty$ bound for $u^+$ in all $\RR^2$. 

In case $K\geq 0$, since $K$ is decreasing, and we are assuming that $u$ 
is a classical solution so that $K$ is continuous, we automatically have 
$K\in L^\infty$. Then, by examining Thm. 2 of Brezis and Merle 
[\BrezisMerle], see also [\ChaLi], we again conclude that $u^+$ 
is uniformly bounded above. 

With $u^+ \in L^\infty$, we now proceed as in the proof 
of Lemma 1 in [\ChaKieCMP], p. 224, to get  
$$
u(x) = u(0) - {1\over 2\pi} 
\int_{\RR^2} \left( \ln |x -y| -\ln |y|\right) K(y)e^{2u(y)} \dy\, .
\eqno\eqlbl\uIntEQpre
$$
Pulling out the contribution $\propto\ln |x|$ from the integral,
noting that 
$$
\ln {|x-y|\over |x||y|} = 
\ln \left |{x\ \over |x|^2 } -{ y\ \over |y|^2}\right|  
\eqno\eqlbl\lnJUGGLE
$$ 
and recalling the definition of the Kelvin transform gives us (\uIntEQ).\qed

\smallskip
\noindent
{\it Proof of Proposition 3.5:} 
Let $|x|\geq 4$. We define, for given $x$, the set 
$$
D_x = \{ y:\, |x|/2 \leq |y|\leq 2|x|\ {\rm and}\ |x-y|\leq 4\}
\eqno\eqlbl\Dxregion
$$
and split $\RR^2$ accordingly into $\RR^2 = D_x \cup D_x^C$, where
$D_x^C$ is the complement of $D_x$ in $\RR^2$. 
Moreover, for $y \in D_x^C$ we use the decomposition
$D_x^C = {E_x}\cup{F_x} \cup {G_x}$
with
$$
{E_x} = \{ y :\, 2 |y|\leq |x|\}\, ,
\eqno\eqlbl\DxCregiona
$$
$$
{F_x} = \{y :\, |y|\geq 2 |x|\}\, ,
\eqno\eqlbl\DxCregionb
$$
$$
{G_x} = \{y :\, |y| \leq 2|x|\leq 4|y|\ {\rm and}\ |x-y|\geq 4\} \, .
\eqno\eqlbl\DxCregionc
$$
Recall (\lnJUGGLE). Let ${\bf I}_\Lambda$ denote 
the indicator function of the set $\Lambda$. 
It is now readily verified that, with positive
generic constants $C$, 
$$
\Bigl| \ln{|X-Y|}\Bigr| = 
	\left| \ln{|x-y|\over |x||y|}\right| 
\leq \left\{ {C \ln |x| + C|\ln |x-y||\ ;  \qquad\qquad\qquad y \in D_x 
\atop 
C + C|\ln |y|| {\bf I}_{{E_x}}
+ C|\ln |x|| {\bf I}_{{F_x}\cup {G_x}}\ ;
\quad y \in D_x^C}
\right.
\eqno\eqlbl\lnESTIMS
$$
In each of these regions the corresponding 
inequality in (\lnESTIMS) follows by an application
of the triangle inequality, paying due attention to the 
a-priori bounds on $x$, $y$, and $x-y$. 

Thus, with positive generic constants $C$,
$$
\eqalignno{
{1\over \ln |x|}
 & \int_{\RR^2} \left| \ln{|x-y|\over |x||y|}\right| |K(y)|e^{2u(y)} \dy
\cr &
\leq {C\over \ln |x|} \int_{|y|\leq 1}  {|\ln|y||} \dy 
&\cr &
\quad + C \int_{1\leq |y|\leq |x|/2}
  {|\ln|y||\over \ln|x|}  |K(y)|e^{2u(y)} \dy 
&\cr &
\quad + C \int_{|y|\geq 2|x|} |K(y)|e^{2u(y)} \dy 
&\cr &
\quad + C 
\int_{|y-x|\leq 4} {| \ln|x-y||\over \ln |x|}  |K(y)|e^{2u(y)} \dy \, .
&\eqlbl\intESTSabc\cr}
$$
The first term  on the right obviously goes to zero as $|x|\to \infty$. 
The second integral on the right goes to zero as $|x|\to \infty$ 
by the dominated convergence theorem, and because $K e^{2u}\in
L^1(\RR^2)$. The third integral on the right
goes to zero as $|x|\to \infty$  because $K e^{2u}\in
L^1(\RR^2)$. For the fourth integral on the right we need to
distinguish two cases, (i) $K\in L^\infty$ and (ii) $K\not\in L^\infty$.
As for case (i), since $u^+\in L^\infty$, we have  
$Ke^{2u}\in L^\infty$, and so
$$
\eqalignno{
{1\over \ln |x|}
 \int_{|y-x|\leq 4} \bigl| \ln|x-y|\bigr|  |K(y)|e^{2u(y)} \dy 
&
\leq {C\over \ln |x|}
\int_{|y-x|\leq 4} \left|
 \ln|x-y|\right|  \dy 
&\cr &
\leq {C \over \ln |x|} \to 0 
\qquad {\rm as} \ |x|\to\infty \, .
&\eqlbl\intESTSc\cr}
$$
As for case (ii), since then $K(x)<0$ for $|x|> R$, we have
$$
-\Delta u(x) = K(x)e^{2u(x)} 
\leq 0 \ \ {\rm for}\ |x|> R\, ,
\eqno\eqlbl\subharmINEQ
$$
whence $u(x)$ is sub-harmonic for $|x|> R$. Thus, for 
$|x_0| \geq  R +1$, we have 
$$
u(x_0) \leq {1\over \pi } \int_{{\rm B}_1(x_0)} u(y) \dy\, .
\eqno\eqlbl\subharmEST
$$
By Jensen's inequality [\HardyLittlewoodPolya],
$$
e^{2u(x_0)} \leq {1\over \pi} \int_{{\rm B}_1(x_0)} e^{2u(y)} \dy \, ,
\eqno\eqlbl\subharmJENSEN
$$
whence
$$
|K(x_0)|e^{2u(x_0)} \leq 
{1\over \pi} \int_{{\rm B}_1(x_0)} |K(x_0)| e^{2u(y)} \dy \, .
\eqno\eqlbl\subharmKEST
$$
Now, by hypothesis, either (\kappastarCone) or (\kappastarCtwo) holds. 
If (\kappastarCone) holds, then $|K(x_0)| \leq C|K(y)|$ for
all $y$ in ${\rm B}_1(x_0)$, whence 
$$
 \int_{{\rm B}_1(x_0)} |K(x_0)| e^{2u(y)} \dy
\leq  
C \int_{{\rm B}_1(x_0)} |K(y)| e^{2u(y)} \dy \leq C\, ,
\eqno\eqlbl\subharmtotCest
$$
where the second estimate holds by (\tCbound). 
It follows  once again that $Ke^{2u}\in L^\infty$,
and so we are back to (\intESTSc). If (\kappastarCtwo) holds,
then, writing $|K| = |K|^{1/p}|K|^{1/q}$ with $p=P$, $1/p+1/q =1$, 
we have, by H\"older's inequality [\HardyLittlewoodPolya],
$$
\eqalignno{
& \int_{|y-x|\leq 4} \bigl| \ln|x-y|\bigr|  |K(y)|e^{2u(y)} \dy 
&\cr &
\qquad \leq 
 \left(
\int_{|y-x|\leq 4} |K(y)|e^{2qu(y)} \dy \right)^{1/q}
\left(
\int_{|y-x|\leq 4} \bigl| \ln|x-y|\bigr|^p  |K(y)| \dy \right)^{1/p}
\, .\qquad
&\eqlbl\intESTSd\cr}
$$
Since $u^+\in L^\infty$, and since (\kappastarCtwo) holds, we now
have
$$
\eqalignno{
{1\over \ln |x|}
& 
	\int_{|y-x|\leq 4} \bigl| \ln|x-y|\bigr|  |K(y)|e^{2u(y)} \dy 
&\cr &
\leq C
 \left(
\int_{|y-x|\leq 4} |K(y)|e^{2u(y)} \dy \right)^{1/q}
\left(
\int_{|y-x|\leq 4} \bigl| \ln|x-y|\bigr|^p  \dy \right)^{1/p}
&\cr &
%\hskip 4truecm 
\to 0  \qquad {\rm as} \ |x|\to\infty \, ,
&\eqlbl\intESTSe\cr}
$$
and this completes the estimates on the third integral in
(\intESTSabc). 

In total, by Lemma 4.1 and our estimates on the last integral in
(\uIntEQ), we conclude that for any $\eps$ 
there exists a $C(\eps)$ and $R(\eps)$ such that 
$$
e^{2u(x)} \leq C |x|^{-(\cK(u)/\pi) +\eps}\qquad {\rm for} \ |x| >
R(\eps)\, .
\eqno\eqlbl\asympDOM
$$
Recalling now the definition of $\kappa_*$, Proposition 3.5 follows.\qed 

\medskip
\noindent
{\bf Lemma 4.2:} {\it 
	Let $u$ be a classical solution of (\PDE) satisfying
	(\uplusLone). Assume (\tCbound) holds, and that 
	$K$ satisfies (\Kraddec). Moreover, if $K$ is unbounded
	below, assume that either (\kappastarCone) or
	(\kappastarCtwo) hold. Finally, if ${\cal K}(u) = \kappa_*(K)$, 
%In case that $S_g$ is a critical least integrally curved surface, 
	let (\intKeulnsq) be satisfied. Then, uniformly in $x$, 
$$
\lim_{|x|\to \infty} 
\left( u(x) - u(0) + {1\over 2\pi} {\cal K}(u)\ln |x|  \right) = 
{1\over 2\pi} \int_{\RR^2} \ln | y| K(y)e^{2u(y)} \dy\, .
\eqno\eqlbl\asympu
$$
}

\medskip
{\it Proof: } By Proposition 3.5, ${\cal K}(u) \geq \kappa_*(K)$. 
If ${\cal K}(u) = \kappa_*(K)$, then (\intKeulnsq) is satisfied, 
by hypothesis, and this implies that
$\int_{\RR^2} \ln |y| K(y)e^{2u(y)}\dy$ exists. 
If ${\cal K}(u) > \kappa_*(K)$, then by (\asympDOM) and the
definition of $ \kappa_*(K)$, the existence of
$\int_{\RR^2} \ln |y| K(y)e^{2u(y)}\dy$ follows once again. 
By inspecting the estimates of the proof of Proposition 3.5, 
we now conclude, once again by dominated convergence, that
$$
\lim_{X\to 0} \int_{\RR^2} \ln | X-Y| K(y)e^{2u(y)}\dy  
= - \int_{\RR^2} \ln |y| K(y)e^{2u(y)}\dy  \, .
\eqno\eqlbl\domconvLIM
$$
Lemma 4.2 follows.~\qed 

\medskip
\chno=5
\equno=0
\noindent
{\bf V. GLOBAL RESULTS} 
\smallskip

With the help of Lemma 4.2, and noting
(\Kraddec), (\tCbound), we now see that the asymptotic behavior 
of $u$ implies that the integral curvature ${\cal K}(u)$ of $S_g\in
{\bf S}_K$ is strictly positive if $K(x) <0$ for $|x|>R$. In addition,  
it follows trivially from the definition of ${\cal K}(u)$
that ${\cal K}(u) \geq 0$ if $K \geq 0$, with equality holding
if and only if $K\equiv 0$. We summarize this as
 
\smallskip
\noindent
{\bf Lemma 5.1:} {\it 
		Let $u$ be a classical solution of (\PDE)
		satisfying (\uplusLone). Assume (\tCbound) holds. In 
		addition assume that $K$ satisfies (\Kraddec). 
		If  ${\cal K}(u) = \kappa_*(K)$, 
		let (\intKeulnsq) be satisfied. Then the 
		integral curvature ${\cal K}(u)$ of $S_g$ is positive, 
$$
\int_{\RR^2} K(x) e^{2{u}(x)}\dx \geq  0 \, ,
\eqno\eqlbl\totradCURVfin
$$
		with `$\, =$' holding iff $K\equiv 0$. }
% \vfill\eject

We will also need an angular average of $u$.  
In the following, we set $r = |x|$, and we identify 
points in $\RR^2$ with points in $\CC$. We define the radial function 
$$
\ol{u}(r) = {1\over 2\pi} \int_0^{2\pi} u\!\left(re^{i\theta}\right)
\dtheta\, ,
\eqno\eqlbl\AVEu
$$
which is well defined for all $r\geq 0$ because
$u$ is a classical solution. Similarly we define
${\ol K}(r)$. Notice that ${\ol K}(|x|) = K(x)$. 

\medskip
\noindent
{\bf Lemma 5.2:} {\it 
		Let $u$ be a classical solution of (\PDE)
		satisfying (\uplusLone), (\tCbound), with $K$ satisfying
		(\Kraddec). If  ${\cal K}(u) = \kappa_*(K)$, 
		let (\intKeulnsq) be satisfied. Let $\ol{u}$ be 
		defined by (\AVEu). Then there exists a positive constant 
		$c(u) < \infty$ such that 
$$
\bigl|u(x) - \ol{u}(|x|)\bigr| \leq c(u)
\eqno\eqlbl\uaveubound
$$
		for all $x$, and $c(u)$ is the smallest such $c$.
}

\medskip
{\it Proof:} For $|x|\leq R$ the statement is trivial, since $u$ is a 
classical solution. For $|x|>R$, the statement follows from Lemma 4.2.\qed

\medskip
\noindent
{\bf Lemma 5.3:} {\it 
		Let $u$ be a classical solution of (\PDE)
		satisfying (\uplusLone), (\tCbound), with $K$ satisfying
		(\Kraddec). Let $\ol{u}$ be defined by (\AVEu). Then we have
$$
\int_{\RR^2} |K(x)| e^{2\ol{u}(|x|)} \dx < \infty\, .
\eqno\eqlbl\totradCURVfin
$$
		If (\intKeulnsq) holds, then we also have
$$
\int_{\RR^2} \bigl(\ln |x|\bigr)^2 |K(x)| e^{2\ol{u}(|x|)} \dx < \infty\, ,
\eqno\eqlbl\intKradeulnsq
$$
}

\medskip
{\it Proof:} By Jensen's inequality,
$$
e^{\ol{u}(r)} \leq {1\over 2\pi} \int_0^{2\pi} 
e^{u\left(re^{i\theta}\right)} \dtheta\, .
\eqno\eqlbl\JensenESTIM
$$
Upon multiplying (\JensenESTIM) by $2\pi r |\ol{K}(r)|$ and then
integrating over $r$, we get
$$
 \int_0^{2\pi} \int_0^\infty |\ol{K}(r)| e^{2\ol{u}(r)} r\dr  \dtheta  
\leq 
\int_{\RR^2} |K(x)| e^{2 u(x)} \dx \, ,
\eqno\eqlbl\totradCURVbound
$$
which now shows that (\totradCURVfin) holds because of
(\tCbound). Similarly, if (\intKeulnsq) holds, then we can
multiply (\JensenESTIM) by $2\pi r (\ln r)^2 |\ol{K}(r)|$ and 
subsequently integrate the result over $r$ to get
$$
\int_0^{2\pi} \int_0^\infty |\ol{K}(r)| e^{2\ol{u}(r)} (\ln r)^2 r\dr \dtheta  
\leq 
\int_{\RR^2} |K(x)| e^{2 u(x)} (\ln |x|)^2 \dx \, ,
\eqno\eqlbl\totradCURVlnbound
$$
which shows that (\intKradeulnsq) now holds because of (\intKeulnsq).\qed

 \vfill\eject
\medskip
\chno=6
\equno=0
\noindent
{\bf VI. THE COMPARISON FUNCTION}
\medskip

In this section we construct a comparison function for 
$u$, a classical solution of (\PDE) satisfying (\uplusLone), 
(\tCbound), with $K$ satisfying (\Kraddec). In case that 
${\cal K}(u) = \kappa_*(K)$, we assume that (\intKeulnsq) is satisfied. 
Recall that $\ol{u}$ is defined by (\AVEu). 

We first introduce a function $g: [0,\infty)\to\RR$, given by
$$
g(r) = r \int_r^\infty |{\ol K}(s)| e^{2\ol{u}(s)} s (\ln s)^2  \ds 
- r\ln r \int_r^\infty |{\ol K}(s)| e^{2\ol{u}(s)} s \ln s \ds\, ,
\eqno\eqlbl\EULEReqSOL
$$
if $r>0$, while $g(0)$ is given by continuous extension to $r=0$. 
Notice that $g$ is well defined for $r\geq 0$, for,
by Lemma 5.3, the integrals are well defined for all
$r\geq 0$, and $r\ln r$ has a removable singularity at $r=0$. 

\medskip
\noindent
{\bf Lemma 6.1:} {\it 
		The function $g$ defined in (\EULEReqSOL)
		is the unique $C^2(\RR^+)$ solution of
		the inhomogeneous Euler equation
$$
 r^2 g^\ppr(r) - rg^\pr (r) + g(r) = |{\ol K}(r)| e^{2\ol{u}(r)} r^3 \ln r\, ,
\eqno\eqlbl\inhomEULEReq
$$
		under the asymptotic condition
$$
g(r) = o(r)\, \qquad as \quad r\to \infty\, . 
\eqno\eqlbl\asympEULERcond
$$
		Furthermore,  $g$ is eventually positive, 
$$
g(r) \geq 0 \ \ \ {\rm if}\ r > 1\, ,
\eqno\eqlbl\EULERgPOS
$$
and $g$ vanishes at $r=0$,
$$
g(0) = 0\, .
\eqno\eqlbl\EULERgNULL
$$
}

\medskip
{\it Proof:} 
Inserting (\EULEReqSOL) into (\inhomEULEReq) one verifies that 
(\EULEReqSOL) is  a particular solution of (\inhomEULEReq). 
Moreover, since $|K| \geq 0$ and $r < s$, when  $r > 1$ we have 
the bounds $0 < (\ln r)(\ln s) < (\ln s)^2$, which imply
$$
0\leq g(r) \leq r \int_r^\infty |{\ol K}(s)| e^{2\ol{u}(s)} s (\ln s)^2  \ds 
\qquad {\rm for}\ r > 1\, .
\eqno\eqlbl\EULERgbounds
$$ 
The first inequality in (\EULERgbounds) states
positivity (\EULERgPOS), and  both together
prove (\asympEULERcond), for clearly
$$
0\leq \lim_{r\to\infty} {g(r)\over r} \leq  
\lim_{r\to\infty}  
\int_r^\infty |{\ol K}(s)| e^{2\ol{u}(s)} s (\ln s)^2  \ds =0\, ,
\eqno\eqlbl\EULERglittleOHr
$$
the last step as a consequence of Lemma 5.3.  
Moreover, since $g(0)$ is defined by $g(0) = \lim_{r\to 0} g(r)$,
(\EULERgNULL) holds because of Lemma 5.3 and $r\ln r \to 0$ for $r\to 0$.  

The general solution of (\inhomEULEReq) is obtained by adding
to this particular solution the general solution of the homogeneous 
problem $Ar + B r\ln r$, with $A,B$ constants.
By (\asympEULERcond), we conclude that $A=B=0$, and thus also 
uniqueness is shown.\qed

Let $\alpha >0$, and define $R(\alpha)$ as the smallest $R > 0$ 
such that $r- \alpha g(r) > e$ for all $r>R$. By 
(\EULERgNULL),  (\asympEULERcond), and by the continuity
of $r\mapsto r - \alpha g(r)$, it follows that a positive $R(\alpha)$
exists, and that $R(\alpha) - \alpha g(R(\alpha)) = e$. 
We now introduce the family of radial functions 
$f_\alpha: \RR^2\backslash {\rm B}_{R(\alpha)}\to\RR$, given by
$$
f_\alpha(x) = \ln (|x| - \alpha g( |x|)) \, .
\eqno\eqlbl\barrier
$$
Clearly, $f_\alpha(x) > 1$ for $|x|> R(\alpha)$, and
$f_\alpha(x) = 1$ for $|x| = R(\alpha)$. We also introduce
$$
\alpha^*(u)  = 2e^{2c(u)} \, ,
\eqno\eqlbl\alphastar
$$
where $c(u)$ is defined in Lemma 5.2.

\medskip
\noindent
{\bf Lemma 6.2:} {\it 
		Given $u$, $\alpha >\alpha^*(u)$, the function 
		$f_\alpha$ defined in (\barrier) satisfies the 
		partial differential inequality
$$
\Delta f_\alpha(x) + 2 K(x) e^{2u(x)}f_\alpha(x) < 0 \, 
\eqno\eqlbl\barrierCRITERION
$$
		for all $x$ satisfying $|x| > \max\{ 1, R(\alpha)\}$. 
}

\medskip
{\it Proof:} In the following, $g^\pr( r) = \partial_r g( r)$, etc. 
Recall that $r=|x|$. 

By explicit calculation we find
$$
\eqalignno{
{\Delta f_\alpha(x)\over f_\alpha(x)} 
&= {\alpha \left(- r^2 g^\ppr(r) + rg^\pr(r) - g(r) \right)
+ \alpha^2 \left( r g(r) g^\ppr(r) - r g^\pr(r)^2 +g(r)g^\pr (r)\right)
	 \over r(r - \alpha g(r))^2\ln(r - \alpha g(r))}
& \cr\cr& 
= -  { \alpha |{\ol K}(r)| e^{2\ol{u}(r)} \over 
	(1 - \alpha g( r)/r)(1+\ln(1 - \alpha g( r)/r)/\ln r)} 
&\cr\cr&\qquad\qquad\qquad
- { \alpha^2  
( g(r) - rg^\pr(r))^2 \over r^2 (r - \alpha g( r))^2\ln(r - \alpha g( r))}
&\cr\cr& 
< - \alpha | K(x)| e^{2\ol{u}(|x|)} \qquad\qquad\qquad\qquad\qquad
{\rm for}\  r > R(\alpha) \, ,
&\eqlbl\firststep\cr}
$$
the last step by the facts that $r>1$ and 
$\alpha g(r) >0$ for $r>1$, and that $r> R(\alpha)$ and 
 $1 -\alpha g /r > 1/r$ for  $r> R(\alpha)$. 
By (\firststep), Lemma 5.2, and $\alpha>\alpha^*(u)$ 
defined in (\alphastar), we now have
$$
\eqalignno{
\Delta f_\alpha(x) + 2 K(x) e^{2u(x)} f_\alpha(x) 
& <
- \left(\alpha |K(x)| e^{2\ol{u}(|x|)} - 
2K(x) e^{2u(x)}\right)f_\alpha(|x|)
&\cr
&\leq 
-  \left( \alpha |K(x)| e^{-2c(u)} -2K(x) \right) e^{2u(x)} f_\alpha(|x|)
&\cr
& \leq 0 
&\eqlbl\fmaxprinc\cr}
$$
for all $x$ satisfying $|x| > \max\{1,R(\alpha)\}$.\qed
\medskip

 \vfill\eject
\medskip
\chno=7
\equno=0
\noindent
{\bf VII. PROOF OF SYMMETRY THEOREMS 3.6 AND 3.7}
\medskip

In the following, we always understand that $S_g\in {\bf S}_K$, that
$u$ is the associated solution of (\PDE), and that (\intKeulnsq) is assumed 
to be satisfied in case that ${\cal K}(u) = \kappa_*(K)$. Moreover,
if $K$ is unbounded below it is also assumed that either
(\kappastarCone) or (\kappastarCtwo) holds.

By Lemma 4.2, $u(x) \to -\infty$ as $|x|\to \infty$. 
Therefore, and since $u$ is a classical solution, $u$ has a global 
maximum, say at $x^*$. Since $K$ satisfies (\Kraddec), if 
$x\mapsto u(x)$ solves (\PDE), then so does $x\mapsto u(\cR(x))$ 
for any $\cR\in SO(2)$. Therefore, after at most a rotation we 
can assume that our solution $u$ has a global maximum at the point
$x^* = (-|x^*|,0)$, with $|x^*| \geq 0$. 

We now introduce the family of straight lines
$$
T_\lambda = \{x\in\RR^2| x_1 = \lambda\}
\eqno\eqlbl\lambdaline
$$ 
and  the half plane `left of $T_\lambda$,' 
$$
\Sigma_\lambda = \{ x : x_1 < \lambda\}\, .
\eqno\eqlbl\halfplane
$$ 
We denote the reflection of $x$ at $T_\lambda$ by
$$
x^{(\lambda)} = (2\lambda -x_1,  x_2)\, .
\eqno\eqlbl\xrefl
$$ 

\medskip
\noindent
{\bf Lemma 7.1:} {\it 
		For $x_1 \leq \lambda \leq 0$, and in particular for 
		$x\in \Sigma_\lambda$ with $\lambda \leq 0$, we have
$$
K(x) \leq K (x^{(\lambda)}) 
\eqno\eqlbl\KlambdaEST
$$
}
\medskip
{\it Proof:} $K$ satisfies (\Kraddec).\qed
\medskip

We next introduce $u_\lambda (x) = u(x^{(\lambda)})$, and also
$$
v_\lambda (x) = u_\lambda (x) - u(x)
\eqno\eqlbl\vlambda
$$
Clearly, $v_\lambda$ is well defined on $\RR^2$. 
 
\medskip
\noindent
{\bf Lemma 7.2:}  {\it 
		For all $\lambda \in \RR $, $v_\lambda$ vanishes on 
		$T_\lambda$ and at infinity, i.e.
$$
\lim_{|x|\to\infty} v_\lambda (x)  =0 
\eqno\eqlbl\asympv
$$
		uniformly in $|x|$.}
 
\medskip
{\it Proof:}  
Notice that on $T_\lambda$ we have $x^{(\lambda)} = x$, whence
$v_\lambda(x)=0 $ for $x\in T_\lambda$. The vanishing of $v_\lambda$
at infinity is a consequence of Lemma 4.2.\qed
\medskip

We next resort to our comparison function $f_\alpha$.
We pick any $\alpha >\alpha^*(u)$ and introduce the function 
$w_\lambda : \Sigma_\lambda \cup T_\lambda \to \RR$, 
defined as 
$$
w_\lambda (x) = \left\{ {
\, v_\lambda (x)/ f_\alpha(x) \qquad |x| \geq R(\alpha)
 \atop 
\ v_\lambda (x)\qquad \qquad\quad |x| \leq R(\alpha)} \right.
\eqno\eqlbl\wlambdaDEF
$$
Notice that $w_\lambda$ is twice continuously differentiable at all
$x$ with $|x| \neq R(\alpha)$, and continuous as function
of $x\in \Sigma_\lambda \cup T_\lambda$, any $\lambda$. 
It vanishes for $|x|\to\infty$ as well as for $x\in T_\lambda$. 
Therefore, if $w_\lambda (x) <0$ for some $x\in \Sigma_\lambda$,
then $w_\lambda$ will have a global negative minimum in $\Sigma_\lambda$.
Our next Lemma will allow us to initialize the moving planes argument,
and also to finalize it. 

\medskip
\noindent
{\bf Lemma 7.3:}  {\it 
		For each $u$ there exists an $R(u)>0$ such that, 
		if $x_*\in \Sigma_\lambda$ is a minimum point for 
		$w_\lambda$, and $w_\lambda(x_*) <0$, then 
		$|x_*| < R(u)$, independently of $\lambda$. }
\medskip

{\it Proof:} We begin by observing that, in the flat case
$K\equiv 0$, then $u = const.$ and $v_\lambda \equiv 0$ for 
all $\lambda$, so that the claim is trivially true. 

In the non-flat case where $K\not\equiv 0$, 
we prove Lemma 7.3 by contradiction. Thus assume that no such
$R(u)$ exists. Then for any $R$ we can find a  $\lambda \leq 0$
such that $|x_*| > R$, where $x_*\in \Sigma_\lambda$ is a minimum 
point for $w_\lambda$ with $w_\lambda(x_*) <0$. 
In particular, we may  choose $R > \max\{ 1, R(\alpha)\}$. 
At such a  minimum point $x_*\in\Sigma_\lambda$ we have
$\nabla w_\lambda(x_*) =0$ and $\Delta w_\lambda(x_*) \geq 0$, 
and of course $w_\lambda(x_*) <0$. 

Now notice first that the reflected function $u_\lambda$ satisfies the PDE
$$
-\Delta {u_\lambda}(x) = K(x^{(\lambda)}) e^{2u_\lambda(x)}\, .
\eqno\eqlbl\PDElambda
$$
Taking the difference between (\PDElambda) and (\PDE) we get
$$
-\Delta {v_\lambda}(x) = 
	K(x^{(\lambda)}) 	e^{2 u_\lambda(x)}
- 	K(x) 		e^{2 u(x)}  \, .
\eqno\eqlbl\PDEdifference
$$
By the mean value theorem there  exists a number
$\psi_\lambda(x)$ between $u(x)$ and $u(x^{(\lambda)})$ such that 
$$
e^{2{u_\lambda}(x)} -  e^{2u(x)}= 2 v_\lambda(x)
e^{2\psi_\lambda(x)}\, .
\eqno\eqlbl\meanvalthm
$$ 
By (\meanvalthm) and Lemma 7.1, we see that  $v_\lambda$ 
satisfies the partial differential inequality
$$
\Delta v_\lambda(x) + 2K(x) e^{2\psi_\lambda(x)}v_\lambda(x)\leq 0 \, ,
\eqno\eqlbl\PDineqV
$$
for all $x\in \Sigma_\lambda$. With the help of 
(\PDineqV) we now easily find that $w_\lambda$ 
satisfies the partial differential inequality
$$
\Delta w_\lambda (x)
+ 2 {\nabla f_\alpha (x) \over f_\alpha(x)}\cdot \nabla w_\lambda(x)
+ \left( {\Delta f_\alpha (x) \over f_\alpha (x)}  
+ 2 K(x)e^{2\psi(x)}\right) w_\lambda (x) \leq 0
\eqno\eqlbl\PDineqW
$$
for all $x\in \Sigma_\lambda$ for which $|x| > \max\{ 1, R(\alpha)\}$.  
Now by assumption $w_\lambda(x_*) <0$, with
$|x_*| > \max\{ 1, R(\alpha)\}$, and since $f_\alpha (x) >1$
for $|x| > \max\{ 1, R(\alpha)\}$, we also have $v_\lambda(x_*) <0$, 
and this means that $u_\lambda(x_*) < u(x_*)$. 
But then $\psi_\lambda (x_*) \leq  u(x_*)$.  Making use of this and 
of $\nabla w_\lambda(x_*) =0$, from (\PDineqW) we now obtain the inequality
$$
\Delta w_\lambda (x_*)
+ \left( {\Delta f_\alpha (x_*)\over f_\alpha (x_*) }
+ 2 K(x_*)e^{2 u(x_*)}\right) w_\lambda (x_*) \leq 0 \, .
\eqno\eqlbl\PDineqWstar
$$
Using now Lemma 6.2, recalling that $\alpha > \alpha^*(u)$,  
in combination with $w_\lambda(x_*)<0$, we see that (\PDineqWstar)
implies that $\Delta w_\lambda (x_*)< 0$. But this is a contradiction
to $\Delta w_\lambda (x_*)\geq 0$.

Hence $w_\lambda$ has no strictly negative minimum outside the disk
$B_{R(u)}$ with $R(u) = \max \{1,R(\alpha^*(u))\}$. 
This concludes the proof of Lemma 7.3.\qed

\medskip
\noindent
{\bf Corollary 7.4:}  {\it 
		For each $u$, when $\lambda < -R(u)$, 
		then $v_\lambda(x)\geq 0$ for $x\in \Sigma_\lambda$.}

\medskip
{\it Proof:} 
Assume $v_\lambda(x_*) < 0$ for some $x_* \in \Sigma_\lambda$, with 
$\lambda < -R(u)$. Then, since $w_\lambda = v_\lambda/f_\alpha$
for all $x\in \Sigma_\lambda$ with $\lambda < -R(u)$,  and since
$f_\alpha > 1$  for all $x\in \Sigma_\lambda$ with $\lambda < -R(u)$,  
we conclude that $w_\lambda(x_*) < 0$ for $x_*\in \Sigma_\lambda$ with 
$\lambda < -R(u)$. But then, since  $w_\lambda\to 0$ as 
$|x|\to \infty$, and $w_\lambda =0$ on $T_\lambda$,  
we see that $w_\lambda$ attains a negative minimum  
for some $x_*\in \Sigma_\lambda$, with $\lambda < -R(u)$.   
This is a contradiction to Lemma 7.3.\qed

Recall the maximum principle (MP) and the Hopf maximum principle (HMP) 
[\GilbargTrudinger]:

\smallskip
\noindent
{\bf MP}: {\it Let 
$\Delta v(x) + \sum_i b_i(x)\partial_{x_i}v(x) + c(x) v(x) \leq 0$
in $\Omega\subset\RR^n$ and $v\geq 0$. 
If $v(\hat{x}) = 0$ for at least one $\hat{x}\in {\rm int}(\Omega)$, 
then $v\equiv 0$ in all of $\Omega$.}

\smallskip
\noindent
{\bf HMP}: {\it Under the same assumptions as in MP, if
$v\not\equiv 0$ in $\Omega$, and $\partial\Omega$
is smooth with $v|_{\partial\Omega} \equiv 0$, then $\partial
v/\partial \nu <0$, where $\partial v/\partial \nu $ is the 
exterior normal derivative on $\partial\Omega$.}

Notice that no sign condition is being imposed on $c(x)$
as the minimum of $v$ is $0$.
\smallskip

We are now ready for the moving lines. The arguments in our 
ensuing proof of Theorems 3.6 and 3.7 are a straightforward 
modification of those in the proof of Theorem P1 in [\ChaKieCMP]. For 
the convenience of the reader we give the complete argument instead
of listing where to modify the arguments of [\ChaKieCMP]. 

\medskip
{\it Proof of Theorems 3.6 and 3.7:} 
By Lemma 7.4, $v_\lambda(x) \geq 0$ for $\lambda < -R(u)$, 
independently of $\lambda$. We now slide the line $T_\lambda$ 
to the right until we reach a critical value $\lambda_0$, which 
is the largest value of $\lambda$ for which $v_\lambda(x) \geq 0$, 
$x\in\Sigma_\lambda$.

\medskip\noindent
{\it Claim A:} $v_\lambda(x) > 0$
for $x\in\Sigma_\lambda$ with $\lambda <\lambda_0$, and
$\partial_{x_1}u >0$ for $x_1 <\lambda_0$.

\medskip\noindent
{\it Claim B:} $\lambda_0 = -|x^*|$.
\medskip

{\it Proof of Claim A}. 
We begin by establishing the first assertion in Claim A. 
Suppose, for $\lambda<\lambda_0$,
that $v_\lambda(x) =0$
at some point $x\in\Sigma_\lambda$. Since $v_\lambda(x)
\geq 0$ for $x\in \Sigma_\lambda$, if $v_\lambda(x) =0$
the minimum of $v_\lambda(x)$ is achieved in  $\Sigma_\lambda$.
Since (\PDineqV) holds, and $v_\lambda(x) \geq 0$, we can apply the
maximum principle and deduce $v_\lambda(x) \equiv 0$ in $\Sigma_\lambda$.
This means for $\lambda = \lambda_0 - \delta$, some $\delta >0$,
that $u(\lambda_0 - 2 \delta, x_2) = u(\lambda_0, x_2)$.
But $v_\lambda(x) \geq 0$, thus $u(x^{(\lambda)}) \geq u(x)$,
which implies $\partial_{x_1}u \geq 0$ for $x_1 \leq \lambda_0$. This
fact together with the fact 
$u(\lambda_0 - 2 \delta, x_2) = u(\lambda_0, x_2)$
yields $\partial_{x_1}u =  0$ for
$\lambda_0 - 2\delta \leq x_1 \leq \lambda_0$.
In particular, $\partial_{x_1}u =  0$ when $x_1 = \lambda_0 - 2\delta$.
By the Hopf maximum principle and the maximum principle we have
$v_\lambda \equiv 0$ iff $\partial_{x_1} {v_\lambda}=0$ on $T_\lambda$.
Now $\partial_{x_1}{v_\lambda} = -2 \partial_{x_1}u$ for $x_1 = \lambda$.
But since $\partial_{x_1}u =0$ when $x_1 = \lambda_0 - 2\delta$,
we see $\partial_{x_1}{v_{\lambda_0 - 2\delta}} =0$ for
$x_1 = \lambda_0 - 2\delta$ or, which is the same,
on $T_{\lambda_0 - 2\delta}$. Now the Hopf maximum principle
says $v_{\lambda_0 - 2\delta} \equiv 0$. We may repeat this
procedure indefinitely and thus deduce that $u$ is independent
of $x_1$. This is a contradiction, and so the first assertion of
(A) is proved.

As for the second assertion of claim A, note that since $v_{\lambda}>0$
in $\Sigma_\lambda$ for $\lambda <\lambda_0$,
and $v_{\lambda} =0$ on $T_\lambda$, 
by the Hopf maximum principle we get
$\partial_{x_1}{v_\lambda} <0$ on $T_\lambda$.
Since for $x_1 = \lambda$ we have $\partial_{x_1}u =
- (1/2)\partial_{x_1}{v_\lambda}$, we also have
$\partial{x_1}(u) >0$ for $x_1 = \lambda$, with $\lambda < \lambda_0$.
So Claim A is proved. 
\eqed

{\it Proof of Claim B}. From the second assertion in Claim A 
we see $u$ is strictly increasing for $x_1 < \lambda_0$.
By a rotation, we had arranged  that the maximum of $u$ 
is at $(-|x^*|,0)$. It follows that $\lambda_0 \leq  -|x^*|$. 
Thus, to prove (B) we need to rule out the case $\lambda_0 < -|x^*|$. 

Assume $\lambda_0 < -|x^*|$. There are two possibilities. 
Either $v_{\lambda_0} \equiv 0$ or $v_{\lambda_0}\not \equiv 0$.

We will  first rule out the case $\lambda_0 < -|x^*|$, and
$v_{\lambda_0}\not \equiv 0$. Indeed, since $v_{\lambda_0}(x) \equiv 0$
for  $x\in T_{\lambda_0}$ and for $|x|\to\infty$ but 
$v_{\lambda_0}\not \equiv 0$  in $\Sigma_{\lambda_0}$, 
and since $v_{\lambda_0}$ satisfies (\PDineqV), by the 
maximum principle we get $v_{\lambda_0} >0$ for $x_1 <\lambda_0$. 
Hence, by the Hopf maximum principle, $\partial_{x_1}{v_{\lambda_0}} <0$ 
when $x_1 =\lambda_0$. On the other hand, by definition of $\lambda_0$,
there exists a sequence of numbers $\lambda_k$, decreasing
to $\lambda_0$, such that $v_{\lambda_k} <0$, whence also
$w_{\lambda_k} <0$,  and $\lambda_0 < \lambda_k <-|x^*|$. Notice
$w_{\lambda_k}$ is well defined for $\lambda_k < -|x^*|$.
Let $x_k$ be a minimum point for $w_{\lambda_k}$.
Then $w_{\lambda_k}(x_k) <0$ and $\nabla w_{\lambda_k}(x_k) =0$.
As Lemma 7.3 implies $|x_k| < R_(u)$, independently of $\lambda$, 
there exists a subsequence $x_{k_j}\to x^*$ such that 
$\nabla w_{\lambda_0}(x^*) =0$ and $w_{\lambda_0}(x^*) \leq 0$ 
for $x^* = (A,B)$, $A \leq \lambda_0$. This is a contradiction. 
Thus our claim is proved in this case.

We will now rule out the case $\lambda_0 < -|x^*|$ and
$v_{\lambda_0} \equiv 0$. Indeed, in that case 
$u(x_{\lambda_0}) = u(x)$ for $x\in \Sigma_{\lambda_0}$.
But $u$ attains its maximum at $(-|x^*|,0)$ and by part
A of our claim, $\partial_{x_1}u >0$ for $x_1 < \lambda_0$.
Since $\lambda_0 < -|x^*|$, it follows $\partial_{x_1}u = 0$ at
$(|x^*|+ 2\lambda_0,0)$, which again is a contradiction.
\eqed

Thus $\lambda_0 = -|x^*|$. Recall that $\lambda_0$ 
is the largest value of $\lambda$ for which $v_\lambda(x) \geq 0$, 
$x\in\Sigma_\lambda$. Hence $v_{-|x^*|}(x) \geq 0$
for $x\in\Sigma_{-|x^*|}$, whence $u_{-|x^*|}(x) \geq u(x)$. 

We may now repeat this argument by sliding the line
$T_\lambda$ in  from $x_1 = \infty$ to get $u_{-|x^*|}(x) \leq u(x)$.
Putting the two inequalities together we conclude that
$u_{-|x^*|}(x) = u(x)$.
This now implies that $u$ is symmetric with respect to
$T_{-|x^*|}$. Moreover, from the arguments involving the 
Hopf maximum principle we see that any solution is also 
decreasing away from $T_{-|x^*|}$. Recall that $(-|x^*|,0)$ 
is the point of global maximum  of $u$.

Finally, we notice that, if $x^*\neq 0$, then, since
$u$ satisfies (\PDE), and since $K$ is radially symmetric about $(0,0)$, 
we conclude that $K$ is a constant. But if $K$ is a constant, 
then if $u(x)$ is a solution of (\PDE), so is $u(x+x^*)$ for any 
fixed $x^*$. Thus, by a simple translation of the origin to $x^*$ 
we can assume that our solution is in fact symmetric with respect to,
and decreasing away from, $T_0$. On the other hand, 
if $K\not\equiv const.$, then $x^* =0$, and again our
solution is symmetric with respect to, and decreasing 
away from, $T_0$. But if $x^* =0$, then we can repeat
our moving line argument with any other than the
$x_1$-direction, thus we come to the conclusion that
$u$ is symmetric about, and decreasing away from, any
straight line through the origin. This now means that $u$
is radially symmetric about and decreasing away from the origin,
modulo a translation in case that $K=constant$.

This completes the proof of Theorems 3.6 and 3.7.\qed

\smallskip
\chno=8
\equno=0
\noindent
{\bf VIII. PROOF OF EXISTENCE THEOREM 2.1}
\smallskip

We begin with the remark that in the special case of identically 
vanishing Gauss curvature our Theorem~2.1 is obviously true. 
Hence, in the rest of this section we assume that the Gauss
curvature is not identically zero.

In the following we will prove a probabilistic theorem which
implies Theorem 2.1 as immediate corollary. Incidentally, the proof 
also provides us with an algorithm for the {\it construction} 
(in principle at least) of nonradial surfaces. We use the 
methods developed in [\KieCPAM] (see also [\clmpCMPcan]), [\KieJllLMP], 
and [\KieSpo]. For applications to Nirenberg's problem, see~[\KiePHYa]. 

We first introduce some probabilistic notation and terminology. 
In the following, $x_1$, $x_2$, ... denote points in $\RR^2$, not
Cartesian components of $x$. 
Let $\NN$ denote the natural numbers. For each $N\in \NN$, we denote 
the probability measures on $\RR^{2N}$ by $P\bigl(\RR^{2N}\bigr)$. 
For $\varrho^{(N)} \in P(\RR^{2N})$, we denote the associated Radon 
measure by $\what\varrho^{(N)}$. A measure 
$\varrho^{(N)} \in P(\RR^{2N})$ is called absolutely continuous
w.r.t. a measure $\varpi^{(N)}\in P(\RR^{2N})$, written 
$\dvarrho^{(N)}\ll \dvarpi^{(N)}$, if there exists a positive
$d\varpi^{(N)}$-integrable function $f(x_1,...,x_N)$, called 
the density of $\varrho^{(N)}$ w.r.t. $\varpi^{(N)}$, such that 
$\dvarrho^{(N)} = f(x_1,...,x_N)\dvarpi^{(N)}$. 
By $P^s(\RR^{2N})$ we denote the exchangeable probabilities, i.e.
the subset of $P(\RR^{2N})$ whose elements are permutation symmetric 
in $x_1,...,x_N$. The $n^{\rm th}$ marginal measure of 
$\varrho^{(N)} \in P^s(\RR^{2N})$, $n<N$, 
is an element of $P^s(\RR^{2n})$, given by
$$
\varrho^{(N)}_n(\!\dx_1...\dx_n ) 
= \int_{\RR^{2N-2n}}  \varrho^{(N)}(\!\dx_1...\dx_n\dx_{n+1}\cdots\dx_N) \, .
\eqno\eqlbl\margRHO
$$ 

By $\Omega\equiv {(\RR^2)}^\NN$ we denote the infinite Cartesian product 
of the  exchangeable $\RR^2$-valued infinite sequences. 
By $P^s(\Omega)$ we denote the permutation symmetric probability 
measures on $\Omega$. The de Finetti-type result of Hewitt and 
Savage [\HewittSavage] states that each $\mu\in P^s(\Omega)$ is uniquely 
presentable as a convex superposition of product measures, 
i.e., for each $\mu \in P^s(\Omega)$ there exists a unique probability 
measure $\nu(\!\dvarrho|\mu)$ on $P(\RR^2)$, such that 
$$
\mu_n (\!\dx_1 ... \dx_n) =
	\int_{P(\RR^2)} \nu(\!\dvarrho|\mu)\, 
		\varrho^{\otimes n}(\!\dx_1 ... \dx_n) \, , 
\qquad n\in \NN\, ,
\eqno\eqlbl\finettimeas
$$
where $\varrho^{\otimes n}(\!\dx_1 ... \dx_n) \equiv 
\varrho(\!\dx_1){\otimes \cdots \otimes} \varrho(\!\dx_n)$, 
and $\mu_n$ denotes the $n^{\rm th}$ marginal measure of $\mu$.
For de Finetti's original work, see [\Finetti]; see also 
[\Dynkin, \Diaconis, \DiaconisFreedman]. 
We remark that (\finettimeas) coincides with the extremal decomposition 
for the convex set $P^s(\Omega)$, an application of the
Krein-Milman theorem. For details, see [\HewittSavage]. 

To $\varrho\in P(\RR^2)$ we assign the energy 
$$
\cE(\varrho) \equiv {1\over 2} \what\varrho^{\otimes 2}(\ln|x-y|)
= {1\over 2} \int_{\RR^2}\int_{\RR^2} \ln|x-y| \varrho(\!\dx)\varrho(\!\dy) 
\, ,
\eqno\eqlbl\Efctnl
$$
whenever the integral on the right exists. We denote by $P_\cE(\RR^2)$
the subset of $P(\RR^2)$ for which $\cE(\varrho)$ exists.  
For $\mu\in P^s(\Omega)$ the mean energy of $\mu$ is  defined as 
$$
e(\mu) = {1\over 2} \what\mu_2(\ln|x-y|)\, ,
\eqno\eqlbl\meanE 
$$
whenever the integral on the right exists. The following 
proposition, proved in [\RuelleBOOK], characterizes the 
subset of $P^s(\Omega)$ for which (\meanE) is well defined.

\smallskip
\noindent
{\bf Proposition 8.1:} {\it The mean energy of $\mu$, (\meanE), is 
well defined for those $\mu$ whose decomposition measure
$\nu(\!\dvarrho|\mu)$ is concentrated on $P_\cE(\RR^2)$, and in
that case given by 
$$
e(\mu) = \int_{P_\cE(\RR^2)} \nu(\!\dvarrho|\mu)\, \cE(\varrho) \, .
\eqno\eqlbl\meanErep
$$
}

Let $\Upsilon:\RR^2 \to \RR^+$ be an $L^\infty$ function, 
$\Upsilon\not\equiv 0$. 
For some entire harmonic function $H$, which may be 
constant, and all $0 < \gamma < 2$, we assume $\Upsilon$ satisfies  
$$
\int_{{\rm B}_1(y)}
\Upsilon(x)e^{2H(x)} |x - y|^{-\gamma} dx \longrightarrow 0
\qquad {\rm as}\ |y|\to \infty\, . 
\eqno\eqlbl\UpsHgamma
$$
Moreover, we assume that for the same harmonic function $H$ 
and some $q>0$, $\Upsilon$ satisfies 
$$
\int_{\RR^2} \Upsilon(x)e^{2H(x)}|x|^q dx <\infty\, ,
\eqno\eqlbl\UpsHq
$$
and define
$$
q^*(\Upsilon,H) = 
\sup\, \Bigl\{ q>0 :
\int_{\RR^2} \Upsilon(x)e^{2H(x)}|x|^q dx <\infty\,  \Bigr\}\, .
\eqno\eqlbl\qSUPstar
$$   
Given such $H$ and $\Upsilon$, we now define the a-priori measure 
$$
\tau(\!\dx) = \Upsilon(x)e^{2 H(x)}\dx\, 
\eqno\eqlbl\aprioriTAU
$$
on $\RR^2$. Since $\Upsilon$ satisfies (\UpsHq), the integral 
$$
M^{(1)} = \int_{\RR^2} \tau(\!\dx)\, 
\eqno\eqlbl\tottaumass
$$
exists and is called the mass of $\tau$. The 
probability measure associated to $\tau$, given by
$$
\mu^{(1)}(\!\dx) = {1\over M^{(1)}}\, \tau(\!\dx)\, ,
\eqno\eqlbl\MUone
$$
is thus clearly absolutely continuous with respect to $\dx$.

For each $\varrho^{(N)}(\!\dx_1...\!\dx_N) \in P\bigl(\RR^{2N}\bigr)$,
its entropy with respect to the probability measure 
 $\mu^{(1)}(\!\dx_1)\otimes ...\otimes\mu^{(1)}(\!\dx_N)\equiv
\mu^{(1) \otimes N}(\!\dx_1 ...\!\dx_N)$ is defined as
$$
\cS^{(N)}\left(\varrho^{(N)}\right) =
 - \int_{\RR^{2N}} 
\ln \left( {d\varrho^{(N)}\quad \over d\mu^{(1)\otimes N}} \right)
 \varrho^{(N)}(\!\dx_1 ... \dx_N)
\eqno\eqlbl\relS
$$
if $\varrho^{(N)}$ is absolutely continuous w.r.t. $d\tau^{\otimes N}$, 
and provided the integral in (\relS) exists. 
In all other cases, $\cS^{(N)}\left(\varrho^{(N)}\right) = -\infty$. 
In particular, if $\mu_n$ is the $n$-th marginal measure of 
a $\mu \in P^s(\Omega)$, then the entropy of $\mu_n$, $n \in\{1,...\}$, 
is given by $\cS^{(n)}(\mu_n)$, with $\cS^{(n)}$ defined as in (\relS)
with $\varrho^{(n)} = \mu_n$. We also define $\cS^{(0)}(\mu_0) = 0$. 

For each $\mu \in P^s(\Omega)$,  the sequence 
$n\mapsto \cS^{(n)}(\mu_n)$ enjoys the following 
useful properties, proofs of which are found in [\RobinsonRuelle], 
(section 2, proof of proposition 1), see also  [\EllisBOOK, \KieCPAM].

\medskip
\noindent
{\it Non-positivity of $\cS^{(n)}(\mu_n)$:} For all $n$,
$$
  \cS^{(n)}(\mu_n)\leq 0\, .
\eqno\eqlbl\Snegativity
$$

\noindent
{\it Monotonic decrease of $\cS^{(n)}(\mu_n)$:} If $n< n^\pr$, then
$$
  \cS^{(n^\pr)}(\mu_{n^\pr})\leq \cS^{(n)}(\mu_n)\, .
\eqno\eqlbl\Sdecrease
$$ 

\noindent
{\it Strong sub-additivity of $\cS^{(n)}(\mu_n)$:}
For $n^\pr,\, n^\ppr \leq n$, and with 
$\cS^{(-m)}(\mu_{-m}) \equiv 0$ for $m>0$,
$$
\eqalignno{
  \cS^{(n)}(\mu_n) \leq \ \cS^{(n^\pr)}(\mu_{n^\pr}) 
& 
	+ \cS^{(n^\ppr)}(\mu_{n^\ppr}) 
&\cr&
	+ \cS^{(n-n^\pr-n^\ppr)}(\mu_{n - n^\pr - n^\ppr})
        - \cS^{(n^\pr + n^\ppr - n)}(\mu_{n^\pr + n^\ppr - n})  \, .
&\eqlbl\Ssubadditivity\cr}
$$

As a consequence of the sub-additivity (\Ssubadditivity) of
$\cS^{(n)}(\mu_n)$, the limit 
$$
  s(\mu) = \lim_{n\to\infty} \ {1\over n} \cS^{(n)}(\mu_n)\, 
\eqno\eqlbl\meanS
$$
exists whenever $\inf_n\, n^{-1} \cS^{(n)}(\mu_n) > -\infty$; 
otherwise $s(\mu) = -\infty$. The quantity $s(\mu)$ given in
(\meanS) is called the mean entropy of $\mu\in P^s(\Omega)$.
The mean entropy is an affine function, see [\RobinsonRuelle].
This entails the following useful representation, proved in [\RobinsonRuelle].
\smallskip
\noindent
{\bf Proposition 8.2:} {\it The mean entropy of $\mu$, (\meanS),
is given by
$$
s(\mu) = \int_{P(\RR^2)} \nu(\!\dvarrho|\mu)\, \cS^{(1)}(\varrho) \, .
\eqno\eqlbl\meanSrep
$$}
\smallskip

Next, identifying each $x_k\in \RR^2$ with the corresponding $z_k\in \CC$, 
we recall the definition of the alternant $\Delta^{(N)}(x_1,...,x_N)$, 
$$
\Delta^{(N)}(x_1,...,x_N)
 =
 \prod_{1\leq i < j \leq N} (z_i - z_j)\, .
\eqno\eqlbl\alternante
$$
Clearly, 
$$
\bigl|\Delta^{(N)}\bigr|(x_1,...,x_N) =
 \prod_{1\leq i < j \leq N} |x_i - x_j| \, .
\eqno\eqlbl\alternanteMOD
$$
We also recall the definition of $q^* >0$ in (\qSUPstar) and define 
$$
\beta^*(\Upsilon,H) = - 2q^*\, .
\eqno\eqlbl\betaSUPstar
$$
For $\beta \in (\beta^*,4)$, and $N\in \NN$,
we now introduce the probability measure $\mu^{(N)}$ on $\RR^{2N}$ by
$$
\mu^{(N)}(\!\dx_1\cdots\dx_N) \equiv {1\over M^{(N)}(\beta)} 
%\prod_{1\leq i < j\leq N} {\left| x_i - x_j \right|}^{-\beta / N} 
\bigl|\Delta^{(N)}\bigr|^{-\beta / N} (x_1,...,x_N)
\prod_{1\leq \ell\leq N} \tau(\!\dx_\ell)\, 
\eqno\eqlbl\canonMU
$$
if $N>1$, and $\mu^{(N)}\equiv \mu^{(1)}$ given in (\MUone) if $N=1$. 
The next Lemma asserts that (\canonMU) is well defined for all $N\in \NN$, 
and all $\beta \in (\beta^*,4)$. 

\smallskip
\noindent
{\bf Lemma 8.3:} {\it For all $\beta\in (\beta^*,4)$, 
the  measure (\canonMU) satisfies $\dmu^{(N)}\ll \dtau^{\otimes N}$. 
Moreover, for the associated density we have
$\dmu^{(N)}/\dtau^{\otimes N}\in L^p(\RR^{2N},\dtau^{\otimes N})$, 
with $p\in [1,\beta^*/\beta)$ when $\beta < 0$, $p\in [1,\infty]$ 
when $\beta = 0$, and $p\in [1,4/\beta)$ when $\beta > 0$. 
}  

\medskip
\noindent
{\it Proof of Lemma 8.3:} 
First, if $\beta =0$, or $N=1$, the claim is obviously true. 

If $N>1$ and $\beta \in (\beta^*,0)$, we make use of the inequality
$$
|x_i - x_j|  \leq (|x_i| +2) (|x_j| + 2) \, ,
\eqno\eqlbl\xiMINUSxjEST
$$
valid for any two $x_i\in \RR^2$ and $x_j\in \RR^2$. Inequality
(\xiMINUSxjEST) is a consequence of the triangle inequality 
$|x_i - x_j| \leq |x_i| + |x_j|$, the fact that $|x|< |x|+2$, and
finally the fact that $r + s   < sr$ when both $r> 2$ and $s> 2$. 
To verify this last inequality, use that $2+r < 2r$ whenever $r>2$, 
so that when $r>2$ and $s>2$, we have $r+s = r + 2 + \eps < 2r + \eps  = 
(2+\eps)r - \eps r +\eps = sr -\eps (r-1) <sr$.
With the help of (\xiMINUSxjEST) we now have for $\beta <0$, 
$$
\eqalignno{
M^{(N)}(\beta) & 
= \int_{\RR^{2N}}
\bigl|\Delta^{(N)}\bigr|^{-\beta / N} (x_1,...,x_N)
%\prod_{1\leq i < j\leq N} {\left| x_i - x_j \right|}^{-\beta / N} 
	\prod_{1\leq k\leq N} \tau(\!\dx_k) 
& \cr &
\leq 
\int_{\RR^{2N}}	
	\prod_{1\leq i \leq N} 	
		{\left(2 +| x_i|\right)}^{-\beta / 2} \tau(\!\dx_i) 
= \left(\int_{\RR^{2}}	
	{\left(2 +|x|\right)}^{-\beta / 2} \tau(\!\dx)\right)^N\, .
&\eqlbl\MminusEST\cr}
$$
The last integral exists, by hypothesis (\UpsHq). This proves  
$\dmu^{(N)}\ll \dtau^{\otimes N}$ for $\beta\in (\beta^*,0)$. 

If $N>1$ and $\beta\in (0,4)$, we use the
inequality between arithmetic and geometric means [\HardyLittlewoodPolya], 
permutation invariance (twice), and H\"older's inequality
[\HardyLittlewoodPolya], to get 
$$
\eqalignno{
M^{(N)}(\beta) & 
= \int_{\RR^{2N}}
\bigl|\Delta^{(N)}\bigr|^{-\beta / N} (x_1,...,x_N)
%\prod_{1\leq i < j\leq N} {\left| x_i - x_j \right|}^{-\beta / N} 
	\prod_{1\leq k\leq N} \tau(\!\dx_k) 
& \cr &
\leq 
\int_{\RR^{2N}}	{1\over N} 
	\sum_{1\leq i \leq N} 	
		\prod_{\scriptstyle 1\leq j\leq N \atop (j\neq i)}
	{\left| x_i - x_j \right|}^{-\beta / 2} 
	\prod_{1\leq k\leq N} \tau(\!\dx_k) 
& \cr &
= \int_{\RR^{2}}\int_{\RR^{2(N-1)}}
	\prod_{2\leq j\leq N }
	{\left| x_1 - x_j \right|}^{-\beta /2 } \tau(\!\dx_j) \tau(\!\dx_1) 
&\cr&
= \int_{\RR^{2}}\left(\int_{\RR^{2}}
{\left| x - y \right|}^{-\beta /2 } 
\tau(\!\dy)\right)^{N-1} \tau(\!\dx) 
&\cr&
\leq \left(
\sup_x 
\int_{\RR^2} {\left| x - y \right|}^{-\beta /2 }\tau(\!\dy)\right)^{N-1}  
M^{(1)} 		% \int_{\RR^{2}} \tau(\!\dx)
 \, .
&\eqlbl\arithmgeomEST\cr}
$$
In the last step we used that for $\beta \in (0,4)$, we have
$$
\int_{\RR^2} |x - y|^{-\beta/2}\tau(\!\dy)
< M^{(1)} + \Psi_\beta(x)\, ,
\eqno\eqlbl\convoluteEST
$$
with
$\Psi_\beta: x\mapsto \int_{B_1(x)} |x -y|^{-\beta/2}\tau(\!\dy)\in C^0(\RR^2)$
(because $\Upsilon\in L^\infty$) and $\Psi_\beta(x)\to 0$ for $|x|\to\infty$
(by hypothesis (\UpsHgamma)). 
This proves  $\dmu^{(N)}\ll \dtau^{\otimes N}$ for $\beta\in (0,4)$. 

By repeating now the same chains of estimates with 
$p\beta$ in place of $\beta$, one concludes that 
$\dmu^{(N)}/\dtau^{\otimes N}\in L^p(\RR^{2N},\dtau^{\otimes N})$ 
for all $p\in [1,4/\beta)$ when $\beta >0$, respectively
all $p\in [1,\beta^*/\beta)$ when $\beta <0$.\qed

\smallskip
We now come to the main theorem of this section. It addresses the
limiting behavior of $\mu_n^{(N)}$  as $N\to\infty$, with $n$
arbitrary but fixed. 

\medskip
\noindent
{\bf Theorem 8.4:} {\it 
		The sequence of probability measures
		$N\mapsto \mu^{(N)}_n(\!\dx_1...\dx_n\, )$ is
		the union of weakly convergent subsequences, in
		the sense that there exist disjoint sequences
		${\bf E}_\ell= \{N_\ell(k)\}_{k\in \NN}$, 
		${\bf E}_\ell \cap {\bf E}_{\ell^\pr}=\emptyset$,
		for $\ell \neq{\ell^\pr}$, such that for each $\ell$, the map
		$k\mapsto \mu^{(N_\ell(k))}_n(\!\dx_1...\dx_n )$ 
		converges weakly in the sense of probability measures, 
		with densities w.r.t. $\dtau^{\otimes n}$ converging 
		weakly in $L^p(\RR^{\, 2n}, \dtau^{\otimes n})$, 
		for all $p\in [1,\infty)$. 

		Let $\mu_n^{\ell}$ denote the weak limit point of 
		such a subsequence. Then there exists a unique 
		$\mu^\ell\in P^s(\Omega)$ (of which $\mu_n^\ell$ 
		is the $n$-th marginal), and $\mu^\ell$ has its
		decomposition measure $\nu\left(\!\dvarrho|\mu^\ell\right)$ 
		concentrated on the subset of 
		$P(\RR^2)\cap \bigcup_{p>1} L^p(\RR^{2},\dtau)$,
		whose elements minimize the functional
$$
\cF_\beta(\varrho) = \beta\cE(\varrho) - \cS^{(1)}\left(\varrho\right) \, .
\eqno\eqlbl\VP
$$ 
}

\bigskip
\noindent
{\bf Remark 8.5:} {\it Notice that Theorem 8.4 asserts that $\cF_\beta$ does 
	have a minimizer $\varrho_\beta^{\phantom{a}}\in P_\cE$. If it 
	can be shown that (\VP) has a unique minimizer, say 
	$\varrho^{\phantom{a}}_{\beta}$, then 
	in fact we have convergence to a product measure, 
$$
\lim_{N\to\infty} \mu^{(N)}_n(\!\dx_1...\dx_n ) = 
{\bigotimes_{1\leq k\leq n}}\,\,
\varrho^{\phantom{a}}_{\beta}(\!\dx_k)\, ,
\eqno\eqlbl\margMUlim
$$
	weakly in 
	$P(\RR^{2n})\cap L^p(\RR^{2n},\dtau^{\otimes n})$ 
	for any $p\in [1,\infty)$. 
} 
\medskip

Before we  prove Theorem 8.4, we show that
Theorem~2.1 is a corollary of Theorem~8.4.

\smallskip
\noindent
{\it Proof of Theorem 2.1:} Assume that all hypotheses of Theorem 8.4
are fulfilled. Then (\VP) has a solution for all $\beta \in (\beta^*,4)$. 
The minimizers of (\VP) are of the form $\varrho(\!\dx) = \rho(x)\dx$,
with $\rho$ satisfying the Euler-Lagrange equation 
$$
\rho(x) = 
{ 			\Upsilon(x)
\exp\left( -{\beta }\int_{\RR^2} \ln|x-y|\, \rho(y)\dy + 2 H(x)\right) 
	\over \int_{\RR^2} \Upsilon(x)
\exp\left(-{\beta } \int_{\RR^2}\ln|x-y|\, \rho(y)\dy + 2 H(x) \right)
		\dx}\, .
\eqno\eqlbl\fixptEQ
$$

Recall that  $\Upsilon \geq 0$, by hypothesis. 
If $\beta \in (0,4)$, we now identify $\Upsilon$ with a (positive) Gauss 
curvature function, $K\equiv \Upsilon$, and if  $\beta \in (\beta^*,0)$ 
we identify $-\Upsilon$ with a (negative) Gauss curvature function, 
$K\equiv -\Upsilon$. In either case, $K$ satisfies the hypotheses of
Theorem~2.1. We also identify $\beta\pi$ with the integral Gauss curvature, 
$$
\kappa = \beta \pi \, ,
\eqno\eqlbl\kabepi
$$
and we notice that $\beta^*\pi = \kappa^*$, defined in (\kappaSUPstar).

We now pick a corresponding solution of (\fixptEQ), say 
$\rho_{_{H,\beta}}$, which exists by Theorem~8.4. With the help of this  
$\rho_{_{H,\beta}}$ we define, for all $x\in \RR^2$, the function 
$$
U_{H,\kappa} (x) = H(x) - 
{\beta\over 2} \int_{\RR^2}\ln|x-y|\, \rho_{_{H,\beta}}(y)\dy  + U_0 \, ,
\eqno\eqlbl\uofrho
$$
the constant $U_0$ being uniquely determined by the requirement that 
$$
\int_{\RR^2} K(x)e^{ 2 U_{H,\kappa}(x)}\dx = \kappa\, .
\eqno\eqlbl\Unull
$$
By Theorem 8.4, $\rho_{_{H,\beta}}\in L^p(\RR^2,\dtau)$ for 
all $p\in [1,\infty)$, whence 
$U_{H,\kappa} \in W_{\rm loc}^{2,p}\cap L_{\rm loc}^\infty$.
With $\Delta \ln|x-y| = 2\pi \delta (x-y)$ it now follows 
that $u(x) = U_{H,\kappa}(x)$ is a distributional solution 
of (\PDE) for the prescribed Gauss curvature function $K$, 
with $K$ satisfying (\KHgamma), (\KHq), and $u$ 
satisfying the asymptotics (\Uasymp). 

For the subset of $K\in C^{0,\alpha}(\RR^2)$ we can
bootstrap to $U_{H,\kappa} \in C^{2,\alpha}(\RR^2)$ by using elliptic
regularity, thus obtaining an entire classical solution of (\PDE).  
For the further subset of $K$ satisfying also (\Krotinv), this
classical solution obviously 
breaks the radial symmetry if $H\not\equiv constant$.
Finally, for the further subset of $K$ satisfying (\Kraddec), 
a straightforward estimate shows that (\KHgamma) is redundant, then.

This concludes the proof of Theorem~2.1.~\qed
\smallskip
We now prepare the proof of Theorem 8.4. 
Let $\Pi(\RR^{2N})$ denote the subset of $P(\RR^{2N})$ whose elements
are absolutely continuous  w.r.t. $\dtau^{\otimes N}$, 
having density $\dvarrho^{(N)}/ \dtau^{\otimes N}
		\in \bigcup_{p>1} L^p(\RR^{2N},\dtau^{\otimes N})$. 
On $\Pi(\RR^{2N})$  we define the functional
$$
\cF_\beta^{(N)}\bigl(\varrho^{(N)}\bigr) = 
	\beta \what\varrho^{(N)}\bigl(
\ln \bigl|\Delta^{(N)}\bigr| \bigr) 
-	N\cS^{(N)}\bigl(\varrho^{(N)}\bigr)\, .
\eqno\eqlbl\FE
$$

\medskip
\noindent
{\bf Lemma 8.6:} {\it For each $\beta\in (\beta^*,4)$, the functional (\FE) 
takes its unique minimum at the probability measure (\canonMU), i.e.,
$$
\min_{\varrho^{(N)} \in \Pi(\RR^{2N})}
 \cF_\beta^{(N)}\bigl(\varrho^{(N)}\bigr) = 
\cF_\beta^{(N)}\bigl(\mu^{(N)}\bigr) \, . 
\eqno\eqlbl\GibbsVP
$$
Moreover,
$$
\cF_\beta^{(N)}\bigl(\mu^{(N)}\bigr) = 
- N \ln \what\mu^{(1)\otimes N} 
\biggl(\, \bigl|\Delta^{(N)}\bigr|^{-\beta/N}\biggr) \, .
\eqno\eqlbl\FEexplicit
$$
For $\beta \geq 4$, and for $\beta < \beta^*$, (\FE) is unbounded below.} 

\smallskip
\noindent
{\it Proof of Lemma 8.6:} 
Since
$\ln \bigl|\Delta^{(N)}\bigr| \in L^p(\RR^{2N},\dtau^{\otimes N})$ for all 
$p\in [1,\infty)$ by Lemma 8.3,
  the integral $\cF_\beta^{(N)}\bigl(\mu^{(N)}\bigr)$ 
is well defined for $\beta \in (\beta^*,4)$. 
The identity (\FEexplicit) is readily verified by explicit 
computation. The Gibbs variational principle (\GibbsVP) in turn is 
just convex duality [\Rockafellar], verified by the standard convexity 
argument (cf., [\EllisBOOK], proof of Proposition I.4.1). Thus, 
rewriting (\FE) as
$$
\cF_\beta^{(N)}\bigl(\varrho^{(N)}\bigr) = \int_{\RR^{2N}} 
\ln \left( {\dvarrho^{(N)}\quad \over \dmu^{(N)}} \right)
 {\dvarrho^{(N)}\quad \over \dmu^{(N)}} \!\dx_1 ... \dx_N
\eqno\eqlbl\relSmuN
$$
and using now $x\ln x \geq x-1$, with equality iff $x=1$, we find that
$\cF_\beta^{(N)}\bigl(\varrho^{(N)}\bigr) \geq 0$, with equality holding
if and only if $\varrho^{(N)} = \mu^{(N)}$. This proves Lemma 8.6 for
$\beta \in (\beta^*, 4)$. 

Now let $\beta \geq 4$, or $\beta <\beta^*$. Assume that
$M^{(N)}(\beta)$ is finite. Then, by (\FEexplicit) 
and by the Gibbs variational principle (\GibbsVP), 
we have
$\min_\varrho \cF_\beta^{(N)}\bigl(\varrho^{(N)}\bigr) 
= - N\ln M^{(N)}(\beta)$.  
However, a simple scaling argument shows that 
$M^{(N)}(\beta \geq 4) > C$ for any $C$, and similarly we have
$M^{(N)}(\beta < \beta^*) >C$ for any $C$, by definition of $\beta^*$. 
This verifies the unboundedness below of (\FE) for $\beta \geq 4$
and $\beta <\beta^*$.\qed

% \vfill\eject

Lemma 8.6 and Lemma 8.3 entail 
\smallskip
\noindent
{\bf Lemma 8.7:} {\it The function $\beta \mapsto F(\beta)$ defined by
$$
F(\beta) \equiv 
\inf_{\varrho \in \Pi(\RR^{2})}\cF_\beta(\varrho)
\eqno\eqlbl\FnullDEF
$$
is continuous for all $\beta\in (\beta^*,4)$. }
% \vfill\eject

\smallskip
\noindent
{\it Proof of Lemma 8.7:} Gibbs' variational principle (\GibbsVP) 
evaluated with a trial product measure 
$\varrho^{(N)} = \varrho^{\otimes N}\in P(\RR^{2N})$, with 
$\varrho \in P(\RR^{2})\cap L^p(\RR^2,\dtau)$ for some $p>1$, gives us
$$
	{1\over N^2} \cF_\beta^{(N)}\left(\mu^{(N)}\right)
% {1\over N} \ln  \left({M}^{(1)}(\beta)^N / {M}^{(N)}(\beta)\right)  
\leq  
	{1\over N^2} \cF_\beta^{(N)}\bigl(\varrho^{\otimes N}\bigr) 
= 
	\Bigl(1-{1\over N}\Bigr) \beta \cE(\varrho ) - \cS^{(1)}(\varrho) 
\eqno\eqlbl\trialPRODmeas
$$
for  all $\varrho \in P(\RR^{2})\cap L^p(\RR^2, \dtau)$, 
$p > 1$, and $N>1$. Now, by (\MminusEST) and (\arithmgeomEST), the left 
side in (\trialPRODmeas) is uniformly bounded below. Letting 
$N\to \infty$ in (\trialPRODmeas) we obtain a lower bound for 
$\cF_\beta(\varrho)$, uniformly over 
$P(\RR^{2})\cap L^p(\RR^2,\dtau)$, $p>1$, for each $\beta \in (\beta^*,4)$. 
Thus,
$$
\eqalignno{
\beta \cE(\varrho) - \cS^{(1)}(\varrho) 
 & \geq \limsup_{N\to \infty}{1\over N^2}\cF_\beta^{(N)}\left(\mu^{(N)}\right)
&\cr &
\geq \liminf_{N\to \infty}{1\over N^2} \cF_\beta^{(N)}\left(\mu^{(N)}\right)
\geq f_0(\beta)\, ,
&\eqlbl\Flowbound\cr}
$$
with
$$
f_0(\beta) = 
\cases{
\displaystyle
	- \ln \int_{\RR^2} (2+ |x|)^{-\beta/2}\mu^{(1)}(\!\dx)\, 
\hskip 1.5truecm &{\rm for }\  $\beta \leq 0\, ,$
	\cr\cr
\displaystyle
	- \ln\, \sup_x\int_{\RR^2} |x - y|^{-\beta/2}\mu^{(1)}(\!\dy)
	\, \qquad & {\rm for }\ $\beta \geq 0\, .$}
\eqno\eqlbl\fnull
$$
Recalling (\VP), this proves that $\cF_\beta$ 
is bounded below for $\beta \in (\beta^*,4)$. 

Having a lower bound, continuity of $F$ now follows from the definition 
of $F$. For assume that $F$ is discontinuous at $\beta_0 \in (\beta^*,4)$.  
Without loss of generality, we can assume $F(\beta_0^-) > F(\beta_0^+)$. 
(The reverse case $F(\beta_0^-) < F(\beta_0^+)$ is treated
essentially verbatim.) 
Now let $\beta = \beta_0 +\eps$. Clearly, for each $\eps$
we can find a minimizing sequence $\{\varrho_k\}_{k\in \NN}$ 
(depending on $\eps$) such that $\cF_{\beta_0+\eps}(\varrho_k) 
< F(\beta^+)+\delta$ if $k>M(\delta)$. Pick a sufficiently small 
$\delta$ and select~a $\varrho_{*} \in \{\varrho_k\}_{k>M(\delta)}$.
Insert this $\varrho_{*}$ into $\cF_{\beta_0 -\eps}$. 
Using $\cF_\beta = \beta\cE\! -\! S^{(1)}$ one gets, 
for any $\eps$ and~$\delta$,
$$
F(\beta_0 -\eps) \leq \cF_{\beta_0 -\eps}(\varrho_{*}) 
= \cF_{\beta_0+\eps}(\varrho_{*}) - 2\eps \cE(\varrho_{*}) 
\leq F(\beta_0 + \eps) +\delta - 2\eps \cE(\varrho_{*}) 
\eqno\eqlbl\FnullCnull
$$
Letting $\eps\to 0$  and $\delta \to 0$ we obtain 
$F(\beta_0^-) \leq F(\beta_0^+)$, which is a contradiction.\qed

\medskip
Taking the infimum over $\varrho$ in (\Flowbound), 
letting $N\to\infty$, and noting Lemma 8.7, gives
\smallskip\noindent
{\bf Proposition 8.8:} {\it For all $\beta\in (\beta^*,4)$, 
$$
\limsup_{N\to \infty}{1\over N^2} \cF_\beta^{(N)}\left(\mu^{(N)}\right)
\leq F (\beta)\, .
\eqno\eqlbl\upbound
$$
}
\smallskip

Proposition 8.8 is complemented by a sharp estimate 
in the opposite direction.

\smallskip
\noindent
{\bf Proposition 8.9:} {\it For all $\beta\in (\beta^*,4)$, 
$$
\liminf_{N\to \infty}{1\over N^2} \cF_\beta^{(N)}\left(\mu^{(N)}\right)
\geq F (\beta) \, .
\eqno\eqlbl\lowbound
$$
}
\smallskip

To prove Proposition 8.9, we need to prove that the sequence of 
the $n^{\rm th}$ marginal measures $\mu_n^{(N)}$ is not 
`leaking at $\infty$' as $N\to\infty$. When $\beta >0$, we 
also need to show that the sequences of 
the densities $\dmu_n^{(N)}/\dtau^{\otimes n}$ of these 
marginal measures are uniformly in $L^p(\RR^{2n},\dtau^{\otimes n})$ 
for $N>N_n(\beta)$. However, since  it gives a-priori regularity, 
we prove uniform $L^p$ bounds for all $\beta\in(\beta^*,4)$. We
remark that when $\Upsilon$ is radial symmetric decreasing, or
has compact support, then many of the following proofs simplify
considerably, some to trivialities. However, since we work with
a minimal set of assumptions on $\Upsilon$, it is unavoidable that
the now ensuing estimates become somewhat more technical.  

We begin by deriving bounds on the expected value of
$\ln\bigl|{\Delta}^{(N)}\bigr|$
w.r.t. $\mu^{(N)}$ which, using permutation symmetry, can be
written in terms of $\mu_2^{(N)}$,
$$
\what\mu^{(N)}\bigl(\ln\bigl|{\Delta}^{(N)}\bigr|\bigr) = 
N(N-1) {1\over 2}\what\mu_2^{(N)}(\ln|x-y|) \, .
\eqno\eqlbl\aveENERGY
$$

\smallskip
\noindent
{\bf Lemma 8.10:} {\it For each $\beta\in (\beta^*,4)$, 
there exist constants $\ul{C}(\beta)$ and $\ol{C}(\beta)$,
independent of $N$, such that for all $N\geq 2$ we have the 
estimates
$$
\ol{C}(\beta)
	\geq 
\beta \what\mu^{(1)\otimes 2}(\ln|x-y|)
	\geq
\beta \what\mu_2^{(N)}(\ln|x-y|) 
	\geq  
\ul{C}(\beta)\, .
\eqno\eqlbl\sandwEST
$$
}
%\vfill\eject
\smallskip
\noindent
{\it Proof of Lemma 8.10:} 
The first inequality in (\sandwEST) is implied by our hypotheses 
(\UpsHgamma) and (\UpsHq) that enter our definitions of $\tau$ (\aprioriTAU)
and $\mu^{(1)}$ (\MUone).  

To obtain the second inequality, we study the functions 
$\beta \mapsto f_N(\beta)$, $N>1$, given by
$$
f_N(\beta) 
= - {2\over N-1}\ln \what\mu^{(1)\otimes N} 
\biggl(\, \bigl|\Delta^{(N)}\bigr|^{-\beta/N}\biggr) \, .
%\biggl( \prod_{1\leq i < j \leq N} |x_i - x_j|^{-\beta/N} \biggr)\, 
\eqno\eqlbl\FperNN
$$
for $\beta\in (\beta^*,4)$. Jensen's inequality [\HardyLittlewoodPolya] 
w.r.t. $\mu^{(1)\otimes N}$ applied in (\FperNN) gives us
$$
f_N(\beta) \leq \beta \what\mu^{(1)\otimes 2}(\ln|x-y|)\, .
\eqno\eqlbl\uppFperNNest
$$ 
On the other hand, $N(N-1)f_N(\beta) = 2\cF_\beta^{(N)}(\mu^{(N)})$. 
Therefore, by Lemma 8.6, (\relSmuN), definition (\FEexplicit),
and the negativity of ${\cal S}^{(N)}$, (\Snegativity), we have
$$
f_N(\beta) 
= \beta \what\mu^{(N)}_2(\ln|x-y|) 
	- {2\over N-1} {\cal S}^{(N)}(\mu^{(N)}) 
\geq \beta \what\mu^{(N)}_2(\ln|x-y|) \, .
\eqno\eqlbl\lowFperNN
$$
The second estimate in (\sandwEST) is proved. 

To prove the third estimate in (\sandwEST), we note that
for any $\beta\in (\beta^*,4)$, there exists a small $\eps >0$
such that $(1+\eps)\beta\in (\beta^*,4)$. By
Jensen's inequality w.r.t. $\mu^{(N)}$,
$$
M^{(N)}((1+\eps)\beta)
 \geq 
M^{(N)}(\beta) \exp\left(-{1\over 2} (N-1)
\eps\beta\what\mu_2^{(N)}(\ln|x-y|)\right)\, .
\eqno\eqlbl\jensentwo
$$
Dividing (\jensentwo) by ${M^{(1)}}^N$, taking the
logarithm and then multiplying by $- 2/(N-1)$ gives
$$
f_N\bigl((1+\eps)\beta\bigr) 
\leq f_N(\beta) +\eps\beta\what\mu_2^{(N)}(\ln|x-y|)\, .
\eqno\eqlbl\energyffbound
$$
Now, $f_N(\beta)$ is bounded above and below independently of $N$, $N>1$,
for  
$$
2 F(\beta) \geq (1-N^{-1})f_{N}(\beta ) \geq   2 f_0(\beta)\, ,
\eqno\eqlbl\uplowFperNNest
$$
$\beta\in (\beta^*,4)$, and since $1-N^{-1}\to 1$. 
The first inequality in (\uplowFperNNest) is Proposition 8.8, the 
second is (\Flowbound).
With the help of (\uplowFperNNest), from (\energyffbound) we now
obtain, for $N>1$, 
$$
\eqalignno{
\beta\what\mu_2^{(N)}(\ln|x-y|) 
& \geq {1\over \eps}\Bigl(
f_N\bigl((1+\eps)\beta\bigr) - f_N(\beta)\Bigr) 
&\cr
&\geq 
{2\over(1-N^{-1}) \eps}\Bigl(
f_0\bigl((1+\eps)\beta\bigr) - F(\beta)\Bigr) 
\geq 
\ul{C}(\beta)
&\eqlbl\lowenergybound\cr}
$$
uniformly in $N$, for all $\beta\in (\beta^*,4)$.~\qed

We next  prove a hybrid bound, which
for $N=1$ reduces to the first inequality in~(\sandwEST). 

\smallskip
\noindent
{\bf Lemma 8.11:} {\it For each $\beta\in (\beta^*,4)$, 
$N\geq 1$, there is an $N$-independent $\widetilde{C}(\beta)$ such that
$$
\beta \what\mu^{(1)}\otimes\what\mu_1^{(N)}(\ln|x-y|) \leq 
\widetilde{C}(\beta)\, .
\eqno\eqlbl\hybridEST
$$
}
\smallskip
\noindent
{\it Proof of Lemma 8.11:} For $\beta = 0$ the statement is obvious. 

For $\beta \in (\beta^*,0)$, we have
$$
\eqalignno{
\beta \what\mu^{(1)}\otimes\what\mu_1^{(N)}(\ln|x-y|) 
&
\leq 
 \what\mu_1^{(N)}\left(\int_{B_1(x)}\beta\ln|x-y|\mu^{(1)}(\!\dy)\right) 
&\cr&
\leq 
 \what\mu_1^{(N)}\left( \widetilde{C}(\beta)\right) = \widetilde{C}(\beta)\, .
&
\eqlbl\hybridESTa
\cr}
$$
The first estimate in (\hybridESTa) is obvious, since
 $\beta \in (\beta^*,0)$. The second estimate follows from the fact that 
$\Psi_{\log}:x\mapsto \int_{B_1(x)}\ln|x-y|\mu^{(1)}(\!\dy)\in C^0(\RR^2)$
(because $\Upsilon\in L^\infty$) with $\Psi_{\log} (x)\to 0$ as
$|x|\to\infty$ (by $\bigl|\ln|x-y|\bigr| < |x-y|^{-\gamma}$ on $B_1(x)$,
 $\gamma\in(0,2)$, followed by (\UpsHgamma)). 

For $\beta \in (0,4)$, we use (\xiMINUSxjEST) to estimate
$$
\what\mu^{(1)}\otimes\what\mu_1^{(N)}(\ln|x-y|) 
\leq 
 \what\mu^{(1)}\bigl(\ln(2+|x|)\bigr) 
+
\what\mu_1^{(N)}\bigl(\ln(2+|y|)\bigr) \, .
\eqno\eqlbl\hybridESTb
$$
By (\UpsHq), 
$$
 \what\mu^{(1)}\bigl(\ln(2+|x|)\bigr) = {C}_1 <\infty \, .
\eqno\eqlbl\hybridESTcA
$$
As to estimating $\what\mu_1^{(N)}\bigl(\ln(2+|y|)\bigr)$, if 
$\beta\in (0,2)$ we can pick $p\in (1,2/\beta)$ and apply H\"older's
inequality w.r.t. $\tau(\dx_1)$ followed by obvious $L^\infty$ 
estimates to get the upper bound 
$\what\mu_1^{(N)}\bigl(\ln(2+|y|)\bigr) \leq 
C(\beta){M^{(N-1)}\bigl(\beta^\pr\bigr)/ M^{(N)}(\beta)}$, 
where $\beta^\pr = (1-N^{-1})\beta$ and where
$$
C(\beta)  =
\left(\int_{\RR^2}\bigl(\ln(2+|y|)\bigr)^{p^*}\tau(\dy)\right)^{1/p^*}
\sup_{N} \sup_{x\in\RR^2}
\left(\int_{\RR^2}|x- y|^{-p\beta^\pr}\tau(\dy)\right)^{1/p}
<\infty \, , 
\eqno\eqlbl\hybridESThoelder
$$
and subsequently estimate the ratio of $M$'s uniformly in $N$ in the
manner done below, but when $\beta\in [2,4)$, H\"older's inequality 
does not lead to $L^\infty$ functions and so this road is then blocked. 
However, noting that for $q\in(0,q^*)$ we have, by (\UpsHq), 
$$
\int_{\RR^2}\exp\bigl( q \ln(2+|y|)\bigr)\tau(\!\dy) = C_2 <\infty\, ,
\eqno\eqlbl\hybridESTcB
$$ 
we can use convex duality [\Rockafellar] for `${\exp}$' to get, 
for any $q\in(0,q^*)$ and all $\beta\in (0,4)$,
$$
\eqalignno{
\what\mu_1^{(N)}\bigl(\ln(2+|y|)\bigr) & - 
{M^{(N-1)}\bigl(\beta^\pr\bigr)\over M^{(N)}(\beta)}
\int_{\RR^2}\exp\bigl( q \ln(2+|y|)\bigr)\tau(\!\dy)
&\cr&
\leq
- {1\over q}\left(1 + \ln q  
+\beta^\pr\what\mu_2^{(N)}(\ln|x-y|)\right) \leq C^*(\beta)\, .
&\eqlbl\hybridESTd\cr}
$$
In (\hybridESTd), $C^*(\beta)$ is independent of $N$, by Lemma 8.10.
Hence, it now remains to 
estimate $M^{(N-1)}\bigl(\beta^\pr\bigr)/M^{(N)}(\beta)$
from above uniformly in $N$, for each $\beta\in (0,4)$. 
To carry out this last step, we  regularize $M^{(N)}$ and prove an
$N$-independent upper bound on the `regularized ratio of $M$'s' 
which is independent of the regularization parameter.
\smallskip

We regularize $\ln |x-y|$ by $-V_\eps(x,y)\equiv \pi^{-2}\eps^{-4}
\int_{B_\eps(x)}\int_{B_\eps(y)}\ln |\xi-\eta|\dxi\deta$. 
Let $\cH_\eps$ denote the Hilbert space 
obtained by completing the $C_0^\infty(\RR^2)$ functions with vanishing 
integral, $\int_{\RR^2}f(x)\dx=0$, w.r.t. the positive definite inner product 
$$
\langle f,f\rangle_{_\eps} \equiv 
N^{-1}\beta \int_{\RR^2}\int_{\RR^2} f(x)  V_\eps(x,y)f(y) \dx\dy \, .
\eqno\eqlbl\GaussESTa
$$
If $B^1 \equiv {B_{1/\sqrt{\pi}}(0)}$ denotes the disk of area 1 centered 
at the origin, and $\delta_y(x)$ is the Dirac measure on $\RR^2$ 
concentrated at $y$, we note that 
$$
\delta^\sharp_{y}(x) \equiv \delta_{y}(x) - \chi_{_{B^1}}(x)\, 
\in \cH_\eps\, . 
\eqno\eqlbl\GaussESTb
$$
Accordingly, 
$$
\delta^\sharp_{(N)}(x) 
\equiv \sum_{k=1}^N \delta^\sharp_{x_k}(x) \, \in \cH_\eps\,  
\eqno\eqlbl\GaussESTc
$$
as well.  We now define 
$$
W_\eps(x) 
\equiv \int_{B^1} V_\eps(x,y)\dy -
{1\over 2}\int_{B^1} \int_{B^1}V_\eps(x,y)\dx\dy 
\eqno\eqlbl\GaussESTd
$$
and write 
$$
\tau(\!\dx) = e^{\beta W_\eps(x)}\wtilde\tau(\!\dx)\, .
\eqno\eqlbl\GaussESTe
$$
Note that, unless $q^* > \beta$, $\wtilde\tau$ does not 
have finite mass, but that does {\it not} cause a problem. We write 
$M^{(N)}_\eps(\beta)$ for $M^{(N)}(\beta)$ with $-\ln|x-y|$ replaced 
by $V_\eps(x,y)$. With (\GaussESTa) to (\GaussESTe), we have the identity
$$
M_\eps^{(N)}(\beta) = 
e^{-{1\over 2} \beta V_\eps(0,0)}
\int_{\RR^{2N}} 
	e^{ {1\over 2}
		\left\langle 
			\delta^\sharp_{(N)},\,
			\delta^\sharp_{(N)} 
		\right\rangle_\eps
	}
\prod_{\ell=1}^N \wtilde\tau(\!\dx_\ell)\, .
\eqno\eqlbl\GaussESTf
$$

We now use Gaussian functional integrals [\GlimmJaffe] to rewrite
(\GaussESTf). Minlos' theorem [\GlimmJaffe] asserts that 
$N^{-1}\beta {V}_\eps(x,y)$ is the covariance `matrix'
of a Gaussian probability measure with mean zero, i.e.,
there exists a Gaussian average $\Gave{\ .\ }$ on a  space
of linear functionals $\Phi$ on  $\cH_\eps$, with 
$\Gave{\phi(x)} =0$ and $\Gave{\phi(x)\phi(y)}= N^{-1}\beta {V}_\eps(x,y)$,
where $\phi(x)$ is shorthand for $\Phi(\delta^\sharp_x)$. Using the
generating function [\GlimmJaffe]
$$
\Gave{ e^{\Phi(f)} } = e^{{1\over 2}\left\langle f,f \right\rangle_\eps}\, 
\eqno\eqlbl\GaussESTg
$$
with $f= \delta_{(N)}^\sharp$ given in (\GaussESTc), then integrating 
over $\RR^{2N}$ w.r.t. $\wtilde\tau^{\otimes N}$, we obtain 
$$
M^{(N)}_\eps(\beta) = 
e^{-{1\over 2} \beta V_\eps(0,0)} 
\Gave{\left( \textstyle{\int_{\RR^2} 
e^{\phi(x)} \tilde\tau(\!\dx)}\right)^N} \, .
\eqno\eqlbl\hybridESTe
$$

Jensen's inequality in the form $\langle F^N \rangle \geq 
\langle F^{N-1} \rangle^{N/(N-1)}$
applied to the right side of (\hybridESTe) now gives,
in terms of the $M_\eps$'s, 
$$
M^{(N)}_\eps(\beta) \geq 
M^{(N-1)}_\eps(\beta^\pr) \left(
M^{(N-1)}_\eps(\beta^\pr)\right)^{1/(N-1)}
\, 
\eqno\eqlbl\hybridESTf
$$
for all $\eps$. Hence, we can now let $\eps \to 0$ and then
$N\to\infty$ to obtain
$$
\eqalignno{
\limsup_{N\to\infty}
{M^{(N-1)}\bigl(\beta^\pr\bigr)\over M^{(N)}(\beta)}
& \leq
\limsup_{N\to\infty}
\left(M^{(N-1)}(\beta^\pr)\right)^{- 1/(N-1)} 
&\cr&
\leq {1\over M^{(1)} } e^{\, \inf_\varrho \cF_\beta(\varrho)}
\leq {1\over M^{(1)} } e^{\,  \cF_\beta(\mu^{(1)}) }
&\cr&
={1\over M^{(1)} }
\exp\left({1\over 2}\beta \what\mu^{(1)\otimes 2}(\ln|x-y|)\right)\, .
&\eqlbl\hybridESTg\cr}
$$
By Lemma~8.10, the r.h.s. of (\hybridESTg) exists and is obviously 
$N$-independent. Combining (\hybridESTcB), (\hybridESTd), and 
 (\hybridESTg),  we have
$$
\what\mu_1^{(N)}\bigl(\ln(2+|x|)\bigr) 
\leq  \wtilde{C}_2(\beta)
\eqno\eqlbl\hybridESTh
$$
independently of $N$. By (\hybridESTb), (\hybridESTcA) and (\hybridESTh), 
and setting $\wtilde{C}(\beta) = \wtilde{C}_1(\beta) +\wtilde{C}_2(\beta)$,
Lemma 8.11 is proved also for $\beta\in (0,4)$.~\qed
\smallskip

We now prepare for uniform $L^p$ bounds. 
\smallskip
\noindent
{\bf Lemma 8.12:} {\it For each $n\in \NN$, $\beta \in(\beta^*,4)$, there
exist $N_n(\beta)\in \NN$ and $C(n,\beta) > 0$, such that
for $N>N_n$ the Radon-Nikodym derivative of $\mu_n^{(N)}$ 
w.r.t. $\tau^{\otimes n}$ is bounded by
$$
{\, \dmu_n^{(N)} \over \dtau^{\otimes n} }(x_1,...,x_n) \leq C(n,\beta) 
\bigl|\Delta^{(n)}\bigr|^{-\beta/N}(x_1,...,x_n) \, .
%\prod_{1\leq i < j\leq n} {\left| x_i - x_j \right|}^{-\beta / N} 
\eqno\eqlbl\LPest
$$
}

\medskip
%\vfill\eject
\noindent
{\it Proof of Lemma 8.12:} When $\beta =0$, this is trivial. 

When $\beta \neq 0$, we begin by writing
$$
{\, \dmu_n^{(N)} \over \dtau^{\otimes n} }(x_1,...,x_n)  = 
 {1\over M^{(N)}(\beta)} \, G(x_1,...,x_n) 
\bigl|\Delta^{(n)}\bigr|^{-\beta/N}(x_1,...,x_n) \, ,
%\prod_{1\leq i < j\leq n} {\left| x_i - x_j \right|}^{-\beta / N} 
\eqno\eqlbl\margDENSITY
$$
where 
$$
G(x_1,...,x_n) = 
\int_{\RR^{2(N-n)}} 
\prod_{\ 1\leq i \leq n < j \leq N}  
{\left| x_i - x_j \right|}^{- \beta / N} 
\prod_{n < k < \ell\leq N} \hskip -.2truecm
{\left| x_k - x_\ell \right|}^{-\beta / N} \tau(\!\dx_j) \, .
\eqno\eqlbl\shortcut
$$
Let $[\![\ .\ ]\!]$ denote integer part. We define
$$
 N_n(\beta)= 
\cases{ 
\displaystyle 
	n\left[\!\!\left[{2\beta^*-\beta\over\beta^*-\beta}\right]\!\!\right]
	\quad 
		& {\rm if}\ $\beta \in (\beta^*,0)\, ,$
	\cr\cr
\displaystyle
	n\left[\!\!\left[{8- \beta \over 4 - \beta}\right]\!\!\right]
	\quad 
		&{\rm if}\ $\beta \in (0,4)\, .$
	\cr}
\eqno\eqlbl\NnullDEFstar
$$
Given $\beta \in (\beta^*,4)$, let $N > N_n(\beta)$. 
Then, by H\"older's inequality, 
$$
\eqalignno{
	G(x_1,...,x_n)  
& \leq 
	\left( \int_{\RR^{2(N-n)}} 
	\prod_{\ 1\leq i \leq n < j \leq N}  
	{\left| x_i - x_j \right|}^{-\beta / 2n} \tau(\!\dx_j) 
	\right)^{2n/N}\times
&\cr&\quad
	\left(\int_{\RR^{2(N-n)}} 
	\prod_{n < i < j \leq N } 
	{\left| x_i - x_j \right|}^{-\beta /( N-2n)} 
	\prod_{n< k\leq N }  \!\!\!\!\!\tau(\!\dx_k)
\right)^{1-2n/N}\, .
&\eqlbl\HoelderEST \cr} 
$$
As for the first factor on the r.h.s. of (\HoelderEST), permutation
symmetry gives 
$$
\int_{\RR^{2(N-n)}}\!\! \prod_{\ 1\leq i \leq n < j \leq N}  
{\left| x_i - x_j \right|}^{-\beta / 2n} \tau(\!\dx_j) =
\left( 
\int_{\RR^2} \prod_{1\leq i \leq n} {\left| x_i - x \right|}^{-\beta / 2n} 
\tau(\!\dx)  \right)^{N-n} \!\!\!\!.
\eqno\eqlbl\permutFIRST
$$
By the arithmetic-geometric mean inequality and permutation invariance, 
we have
$$
\eqalignno{
& \int_{\RR^2} \prod_{1\leq i \leq n} {\left| x_i - x \right|}^{-\beta / 2n} 
\tau(\!\dx) 
\leq {1\over n} \sum_{1\leq i \leq n} \int_{\RR^2} 
	{\left| x_i - x \right|}^{-\beta / 2} \tau(\!\dx)  \, .
&\eqlbl\firstnewEST \cr} 
$$
For the r.h.s. of (\firstnewEST), we have the estimates
$$
{1\over n} \sum_{1\leq i \leq n} \int_{\RR^2} 
	{\left| x_i - x \right|}^{-\beta / 2} \tau(\!\dx)  
\leq \cases{\displaystyle
	C_n \int_{\RR^2} {\left(2+ |x|\right)}^{-\beta / 2} \tau(\!\dx)
 &{\rm if} $\beta < 0\, ,$ 
\cr\cr
\displaystyle
	\sup_y \int_{\RR^2} {\left| y - x \right|}^{-\beta / 2} \tau(\!\dx)
 &{\rm if} $\beta \geq 0\, ,$ 
\cr}
\eqno\eqlbl\Best
$$
where $C_n = \max_{i\in\{1,...,n\}} (2+|x_i|)^{|\beta|/2}$. 
By (\permutFIRST), (\firstnewEST), and (\Best)
the first term on the right-hand side 
of (\HoelderEST) is bounded by the $2n(1-n/N)$-th power of the 
right-hand side of (\Best), whence uniformly w.r.t. $N$.

As for the second factor on the r.h.s. of (\HoelderEST), we split off
the $(-2n/N)$-th power. We set $\alpha(N) = (N-n)/(N-2n)$. 
Since $N>N_n$, we have $1 < \alpha(N)< 4/\beta$ if $\beta > 0$ 
and $1 < \alpha(N)< \beta^*/\beta$ if $\beta <0$. 
We also have $\alpha(N)\to 1$ as $N\to\infty$. 
Proceeding as in the proof of Lemma 8.7, we find that
$$ 
\eqalignno{
&\limsup_{N \to \infty }\ \left(
\int_{\RR^{2N-2n}} 
	\prod_{n< i < j \leq N} \!\!
	{\left| x_i - x_j \right|}^{-\beta /( N-2n)} 
	\prod_{n< k\leq N } \!\! \tau(\!\dx_k) \right)^{- 2n/ N} 
&\cr\cr &
\hskip 1truecm = 
\limsup_{N \to \infty }\
\left( M^{(N-n)}\bigl(\alpha(N)\beta\bigr) \right)^{- 2n/ N} 
\quad \leq \quad \left(e^{ F(\beta)}/M^{(1)}\right)^{2n} \, ,
&\eqlbl\secondnewEST \cr} 
$$
which implies an $N$-independent bound on 
$\left( M^{(N-n)}\bigl(\alpha(N)\beta\bigr) \right)^{- 2n/ N}$.
Feeding (\permutFIRST), (\firstnewEST), (\Best), and 
(\secondnewEST) back into (\HoelderEST) we see that
$G(x_1,...,x_n)\leq C M^{(N-n)}\bigl(\alpha(N)\beta\bigr)$. 
This already proves that the density (\margDENSITY) is eventually (if
$N$ is big enough) in any $L^p(\RR^{2n}\dtau^{\otimes n})$, $p<\infty$. 
To prove that 
$\dmu_n^{(N)} / \dtau^{\otimes n} \in L^p(\RR^{2n}\dtau^{\otimes n})$
uniformly in $N$, it remains to estimate the ratio 
$M^{(N-n)}\bigl(\alpha(N)\beta\bigr)/M^{(N)}(\beta)$ 
from above, independently of $N$. 

We  once again can
apply the Gaussian functional integral method used in the proof of
Lemma 8.11. Since $\alpha(N)\beta$ occurs
in the argument of $M^{(N-n)}$ instead of $(1-nN^{-1})\beta$, 
beside Jensen's inequality (now pulling a power $N/(N-n)$ out of the
average), we now also need  a `change of covariance formula,'
see [\GlimmJaffe]. However, having proved
Lemmata 8.10 and 8.11 already, a more direct way is the following. 

Using Jensen's inequality twice in a self-explanatory way, we obtain
$$
\eqalignno{
\hskip -1.5truecm
{M^{(N-n)}\Bigl(\alpha(N)\beta\Bigr)\over M^{(N)}(\beta)}
& \leq {1\over {M^{(1)}}^n }
\exp\left( {n(n-1)\over N}{1\over 2}\beta\, 
		\what\mu^{(1)\otimes 2}(\ln|x-y|)\right) \times
&\cr&\hskip 1truecm
\exp\left(n\left(1 -{n\over N}\right)\beta \,
\what\mu^{(1)}\otimes \what\mu_1^{(N-n),\alpha}(\ln|x-y|)\right) \times
&\cr&\hskip 1.5truecm 
\exp\left(- n\left(1-{n+1\over N}\right)
\alpha(N)\beta\, \what\mu_2^{(N-n),\alpha}(\ln|x-y|)\right)\, ,
&\eqlbl\thirdnewEST \cr} 
$$
where $\what\mu^{(N-n),\alpha}$ stands for (\canonMU) 
with $\alpha(N)\beta$ in place of $\beta$. 
The first exponential factor on the r.h.s. of (\thirdnewEST) is 
bounded above uniformly in $N$ because
$\beta\what\mu^{(1)\otimes 2}(\ln|x-y|)\leq \ol{C}(\beta)$
independently of $N$, by the first inequality in Lemma 8.10; as for the
second exponential factor on the r.h.s. of (\thirdnewEST), 
by re-identifying $N\to N-n$ and $\beta\to\alpha\beta$, 
the $N$-independent upper bound in Lemma 8.11 gives
$\beta\what\mu^{(1)}\otimes \what\mu_1^{(N-n),\alpha}(\ln|x-y|) \leq
\widetilde{C}(\alpha(N)\beta)$.  Since $\alpha(N)\to 1$ as $N\to\infty$, 
the second exponential factor on the r.h.s. of (\thirdnewEST) is 
bounded above uniformly in $N$. 
As for the third exponential factor on the r.h.s. of (\thirdnewEST), 
since $\beta^*< \alpha(N)\beta < 4$, and since (\sandwEST) holds for
all $\beta \in (\beta^*,4)$, by Lemma 8.10 we have 
$\beta \what\mu_2^{(N-n),\alpha}(\ln|x-y|) \geq \ul{C}(\alpha(N)\beta)$. 
Again since $\alpha(N)\to 1$, we now see that also the
third exponential factor on the r.h.s. of (\thirdnewEST) is
bounded above uniformly in $N$. This proves Lemma 8.12.\qed

Lemma~8.12 establishes that for each triple
$n\in\NN,\ \beta\in(\beta^*,4),\ p\in [1,\infty)$ there exists a
$\widetilde{N}_n(\beta,p)$ ( $>N_n(\beta)$) such that 
$\dmu_n^{(N)}/\dtau^{\otimes n}\in L^p(\RR^{2n}, \dtau^{\otimes n})$
uniformly in $N$ when $N>\widetilde{N}_n(\beta,p)$.
Hence, the sequence $N\mapsto\mu_n^{(N)}$ 
is $L^p(\RR^{2n}, \dtau^{\otimes n})$-weakly compact
when $N > \widetilde{N}_n(\beta,p)$, for each $p\in [1,\infty)$.

However, a weak $L^p$ limit point of $\mu_n^{(N)}$ need not be a 
probability measure. Since $\RR^2$ is unbounded, some partial 
mass of the marginals $\mu_n^{(N)}$ of (\canonMU) could 
escape to infinity when $N\to\infty$. We now show 
that this  does not happen by proving tightness of the sequences. 
Recall [\BillingsleyBOOK] that the sequence of probability 
measures $\mu_n^{(N)}$ is {\it tight} if for each $\eps\ll 1$ there
exists $R(\eps)$ such that $\mu_n^{(N)}({\rm B}_{R(\eps)}^{\ n}) > 1-\eps$,
independent of $N$, where ${\rm B}_R^{\ n} \subset \RR^{2n}$ is 
the $n$-fold Cartesian product of the ball ${\rm B}_R\subset\RR^2$ 
that is centered at the origin, having radius $R$.

\medskip\noindent
{\bf Lemma 8.13:} 
{\it For each $n$, the sequence $\bigl\{\mu_n^{(N)}\bigr\}_{N\geq n}$ 
given by (\canonMU) is tight.}

\medskip\noindent
{\it Proof of Lemma 8.13:} Since our marginal measures are 
permutation symmetric and consistent, in the sense that 
$\mu_n^{(N)}(\!\dx^n) = 
\mu_{m}^{(N)}\bigl(\!\dx^n\otimes\RR^{2(m-n)}\bigr)$ 
for $m>n$, it suffices to prove tightness for $n=1$. 

It follows from the definition of $\mu^{(1)}$ that the map
$y\mapsto h(y) \equiv \int_{\RR^2}\ln|y-x|\mu^{(1)}(\!\dx) +C$ 
is continuous and  independent of $N$. The constant
$C$ is chosen so that $h(y)>0$. Moreover, we have
$h(y) \to \infty$ as $|y|\to\infty$, uniformly in $y$. 
Therefore, and by Lemma 8.11, for each positive $\eps \ll 1$, we 
can find $R(\eps)$, independent of $N$, such that for all $N$,
$$
\inf_{x_1\ \not\in \ {\rm B}_{R(\eps)}} h(x_1) 
\geq 
{1\over \eps}\what\mu_1^{(N)}(h(x_1)) \, .
\eqno\eqlbl\Uratio
$$
Let ${\bf I}_\Lambda$ denote the indicator function of the set
$\Lambda$. We then have the chain of estimates
$$
\eqalignno{
\what\mu_1^{(N)}(h(x_1)) 
& \geq 
\what\mu_1^{(N)}\bigl(h(x_1){\bf I}_{\RR^2\backslash {\rm B}_{R(\eps)}}\bigr) 
&\cr & 
\geq
{1\over \eps}\what\mu_1^{(N)}\bigl(h(x_1)\bigr) 
\what\mu_1^{(N)}\bigl({\bf I}_{\RR^2\backslash {\rm B}_{R(\eps)}}\bigr) 
 &\cr &
= {1\over \eps}
\what\mu_1^{(N)}(h(x_1)) 
 \(1 - \mu_1^{(N)}\bigl({\rm B}_{R(\eps)} \bigr) \)\, .
&\eqlbl\chainEST \cr}
$$
Dividing (\chainEST) by $\eps^{-1} \what\mu_1^{(N)}(h(x_1))$ and 
resorting terms slightly gives us
$$
\mu_1^{(N)}\bigl({\rm B}_{R(\eps)} \bigr)  \geq 1 - \eps\, ,
\eqno\eqlbl\tightness
$$
independent of $N$. The proof is complete.\qed

% \vfill\eject

To prove Proposition 8.9 we also need a lower
bound on the mean entropy.  

\smallskip
\noindent
{\bf Lemma 8.14:} {\it For each $\beta\in (\beta^*,4)$, there
exists a $C(\beta)$, independent of $N$,  such that
$$
{1\over N} \cS^{(N)}\left( \mu^{(N)}\right) \geq  C(\beta)\, .
\eqno\eqlbl\SperNlowEST
$$}
\smallskip

\noindent
{\it Proof of Lemma 8.14:} 
By definition (\FE) of $\cF_\beta^{(N)}\bigl(\varrho^{(N)}\bigr)$, 
$$
{1\over N } \cS^{(N)}\left( \mu^{(N)}\right) = 
\beta {1\over N^2 } \what\mu^{(N)}\left(\ln\bigl|{\Delta}^{(N)}\bigr|\right) 
	-  {1\over N^2 } \cF_\beta^{(N)}\left( \mu^{(N)}\right)  \, .
\eqno\eqlbl\SperN
$$
The bound (\SperNlowEST) now follows
from Proposition 8.8, (\aveENERGY) and Lemma 8.10.\qed

\smallskip
\noindent
{\it Proof of Proposition 8.9:} By Lemma 8.13, the sequence of probability 
measures $\{\mu_n^{(N)}\vert N = n,n+1,...\}$ is tight in
$P(\RR^{2n})$ for all $n$. Therefore [\BillingsleyBOOK] we 
can select a subsequence $k\mapsto N_\ell(k)\in \NN$, $k\in \NN$, 
such that for each $n\in \NN$, 
$\mu_n^{(N_\ell(k))}\rightharpoonup\mu_n^\ell\in P(\RR^{2n})$, 
as $k\to\infty$. Since the marginals are consistent
(in the sense defined above in the proof of tightness), 
by Kolmogorov's existence theorem (see [\BillingsleyBOOK] p.228 ff.,
see also [\EllisBOOK] p.301 ff.) 
the infinite family of marginals $\{\mu_n^\ell\}_{n\in \NN}$
now defines a unique $\mu^\ell\in P^s(\Omega)$, 
Furthermore, for $\beta^*<\beta<4$, we have as corollary of
Lemma 8.12 that, for any $n$ and any $p\in [1,\infty)$, 
the sequence $\{\mu^{(N)}_n\vert N = n,n+1,...\}$ is eventually in 
a ball $\{g:\|g\|_{L^p(\RR^{2n})}\leq T \}$, where $T(n,\beta,p)$ is
independent of $N$. Therefore, as  ${k\to\infty}$, after
at most selecting a sub-subsequence (also denoted by 
$k\mapsto N_\ell(k)\in \NN$, $k\in \NN$), we have 
that $\dmu_n^{(N_\ell(k))}/\dtau^{\otimes n}\rightharpoonup
\dmu_n^\ell/\dtau^{\otimes n}$, weakly in 
$L^p(\RR^{2n},\dtau^{\otimes n})$, any $p\in [1,\infty)$.  

We first study convergence of energy. By (\aveENERGY) we have
$$
{1\over N_\ell(k)^2}\what\mu^{(N_\ell(k))}
\bigl(\ln\bigl|{\Delta}^{(N_\ell(k))}\bigr|\bigr)
= \(1 -{1\over N_\ell(k)}\) {1\over 2} \what\mu_2^{(N_\ell(k))}(\ln|x-y|)\, .
\eqno\eqlbl\aveENERGYperNN
$$
Since $\ln |x - y| \in {L^q}(\RR^4,\dtau^{\otimes 2})$,
${1\over q} + {1\over p} =1$, by 
weak $L^p(\RR^4,\dtau^{\otimes 2})$ convergence 
of $\mu_2^{(N_\ell(k))}$, 
$$
{1\over 2} \what\mu_2^{(N_\ell(k))}(\ln|x-y|) \to
{1\over 2} \what\mu_2^\ell(\ln|x-y|) = e(\mu^\ell)\, .  
\eqno\eqlbl\Econverge
$$
Since furthermore $1- {N_\ell(k)}^{-1} \to 1$ as $k\to\infty$, we have
$$
\lim_{k\to\infty}{1\over N_\ell(k)^2}
\what\mu^{(N_\ell(k))}\bigl(\ln\bigl|{\Delta}^{(N_\ell(k))}\bigr|\bigr)
= e(\mu^\ell)\, .  
\eqno\eqlbl\EperNNtomeanE
$$

We now turn to the entropy. We define
$m= N_\ell(k)-[\![N_\ell(k)/n]\!]n$. By subadditivity (\Ssubadditivity) 
and negativity (\Snegativity) we have, for any $n< N_\ell(k)$,
$$
\eqalignno{
 {1\over N_\ell(k)} \cS^{(N_\ell(k))} \left(\mu^{(N_\ell(k))}\right) 
& 
\leq{1\over N_\ell(k)} \left[\!\!\left[ {N_\ell(k)\over n}\right]\!\!\right]
\cS^{(n)} \left(\mu_n^{(N_\ell(k))}\right)
+{1\over N_\ell(k)} \cS^{(m)}
\left(\mu_{m}^{(N_\ell(k))}\right)
&\cr &
\leq {1\over N_\ell(k)}
\left[\!\!\left[{N_\ell(k)\over n}\right]\!\!\right]
\cS^{(n)}\left(\mu_n^{(N_\ell(k))}\right) \, .
&\eqlbl\SsubaddEST\cr }
$$
Clearly, 
$N_\ell(k)^{-1}[\![ N_\ell(k)/n ]\!]\to n^{-1}$. Moreover, for each $n$, 
weak upper semicontinuity of $\cS^{(n)}$ (see [\RuelleBOOK]) gives us
$$
\limsup_{k\to\infty}\ \cS^{(n)}\Bigl(\mu_n^{(N_\ell(k))}\Bigr)  
\leq \cS^{(n)}(\mu_n^\ell)\, .
$$
Therefore, for all $n$, 
$$
\limsup_{k\to\infty} 
{1\over N_k}\cS^{(N_\ell(k))}\Bigl(\mu^{(N_\ell(k))}\Bigr) 
\leq {1\over n} \cS^{(n)}\left(\mu_n^\ell\right)\, .
\eqno\eqlbl\SperNuppEST
$$
Recalling (\meanS) and Lemma 8.14, we see that $s(\mu)$ exists. 
Hence, $n\to\infty$ in (\SperNuppEST) gives
$$
   \limsup_{k\to\infty} {1\over N_\ell(k)}\cS^{(N_\ell(k))}
\left(\mu^{(N_\ell(k))}\right) \leq s(\mu^\ell)
\eqno\eqlbl\SperNmeanS
$$
for each convergent subsequence $\mu^{(N_\ell(k))}\rightharpoonup \mu^\ell$. 

Pulling the estimates (\EperNNtomeanE) and 
(\SperNmeanS) together we find, for any $\beta\in (\beta^*,4)$,
$$
   \liminf_{k\to\infty} {1\over N^2_\ell(k)}\cF_\beta^{(N_\ell(k))}
\left(\mu^{(N_\ell(k))}\right) \geq \beta e(\mu^\ell) - s(\mu^\ell)\, .
\eqno\eqlbl\FperNNaboveF
$$
Recalling now Propositions 8.1 and 8.2, and finally using Lemma 8.7, we have
$$
\beta  e(\mu^\ell) - s(\mu^\ell) =
\int_{P(\RR^2)} \nu(\!\dvarrho|\mu^\ell)\, \cF_\beta(\varrho) \geq 
 F(\beta)\, .
\eqno\eqlbl\FElowBound
$$
By (\FElowBound) and (\FperNNaboveF),
the proof of Proposition 8.9 is complete.\qed

 We remark that, when $\beta <0$, Proposition~8.9 can be proved 
without $L^p$ estimates. Indeed, when $\beta <0$, then 
(\FperNNaboveF) follows already with (\EperNNtomeanE) replaced by
$$
\limsup_{k\to\infty} {1\over 2} \what\mu_2^{(N_\ell(k))}(\ln|x-y|) \leq 
{1\over 2} \what\mu_2^\ell(\ln|x-y|) = e(\mu^\ell)\, ,
\eqno\eqlbl\EperNNabovemeanE
$$
which holds by the weak upper semi-continuity of $\ln |x-y|$ 
and the weak convergence of $\mu_2^{(N)}$ in the sense of measures,
see [\KieCPAM, \KieSpo]. Also the entropy estimate 
in the proof of Proposition 8.9 holds by just such
weak convergence of $\mu_2^{(N)}$.
However, without $L^p$ estimates one loses the a-priori information 
on the regularity of the solutions of (\fixptEQ).

% \vfill\eject
We now prove our main existence theorem. 
\smallskip
\noindent
{\it Proof of Theorem 8.4.}

Combining Propositions 8.8 and 8.9, we conclude that
$$
\lim_{N\to\infty} {1\over N^2} \cF_\beta^{(N)}\bigl(\mu^{(N)}\bigr) = F(\beta)\, .
\eqno\eqlbl\FElim
$$
Recalling (\lowbound) and (\FElowBound), we see that (\FElim)
implies  
$$
\int_{P(\RR^2)} \nu(\!\dvarrho|\mu^\ell)\, \cF_\beta(\varrho) = F(\beta)
\eqno\eqlbl\FEdefinetti
$$
for every limit point $\mu^\ell$ of $\mu^{(N)}$. 
Equation (\FEdefinetti) in turn implies that
the decomposition measure $\nu(\!\dvarrho|\mu^\ell)$ is concentrated 
on the minimizers of $\cF_\beta(\varrho)$; for assume not, then by Lemma
8.7, we would have
$$
\int_{P(\RR^2)} \nu(\!\dvarrho|\mu^\ell)\, \cF_\beta(\varrho) > F(\beta) \, ,
$$
which contradicts (\FEdefinetti). The proof of Theorem 8.4 is complete.\qed

We now are also in the position to vindicate our Remark 8.5. 
By the tightness and weak $L^p$ compactness, the sequence 
$\{ \mu^{(N)}, N=1,2,\ldots \}$ is a union of weakly convergent 
subsequences in $L^p$. If the minimizer $\varrho_\beta$ is 
unique, the set of limit points of $\{ \mu^{(N)}, N\in\NN\}$ 
consists of the single product measure 
$\mu = \varrho_\beta^{\otimes \NN}$.~\eqed

\smallskip
\chno=9
\equno=0
\noindent
{\bf IX. PROOF OF UNIQUENESS THEOREM 2.2 FOR $K\leq 0$}
\smallskip

We conclude this paper with a proof of Theorem~2.2. We do this
by proving the dual version, i.e. uniqueness of solutions of 
(\fixptEQ) when $\beta <0$. 

\smallskip
\noindent
{\bf Theorem 9.1:} {\it 
	For $\beta <0$ the solution $\rho_{\beta,H}$ 
	of the fixed point equation (\fixptEQ) is unique.}

\smallskip
\noindent
{\it Proof:} We introduce operator notation for (\fixptEQ), thus
$$
\rho = {\cP}(\rho)
\eqno\eqlbl\fixptOPeq
$$ 
where $\cP$ indicates that the right side is a probability density
over $\RR^2$. 
Now assume that for given $\beta <0$ and $H$ entire harmonic, 
two solutions of (\fixptOPeq) exist, say $\rho^{(1)}$ and $\rho^{(2)}$. 
Then $\rho_{2,1}\equiv \rho^{(2)} - \rho^{(1)}\in H_0^{-1}(\RR^2)$. 
In particular, $\int_{\RR^2}\rho_{2,1}\dx =0$, and 
$$
-\int_{\RR^2} \int_{\RR^2} 
\rho_{2,1}(x)\ln|x-y|\rho_{2,1}(y)\dx\dy \geq 0
\eqno\eqlbl\posdefness
$$
with equality holding iff $\rho_{2,1}\equiv 0$, cf. [\SaffTotikBOOK]. 

For $\lambda \in [0,1]$, we now define the interpolation density 
$
\rho_\lambda = \rho^{(1)} +\lambda \rho_{2,1}\, .
$ 
Expected value w.r.t. $\cP(\rho_\lambda)$ is denoted by
$$
\langle f\, \rangle(\lambda) = \int_{\RR^2} f(x) \cP(\rho_\lambda)(x)\dx\, .
\eqno\eqlbl\lambdaMEAN
$$ 
We now use (\fixptEQ) for one of the $\rho_{2,1}$ in the l.h.s. of
(\posdefness) and, with the abbreviation
$$
 U_{2,1}(x) =   \int_{\RR^2} \ln|x-y|\rho_{2,1}(y)\dy\, ,
\eqno\eqlbl\UonetwoDEF
$$
find that
$$
\eqalignno{
 -\int_{\RR^2} \int_{\RR^2} 
 \rho_{2,1}(x)\ln|x-y|\rho_{2,1}(y)\dx\dy 
& = - \int_{\RR^2}  \int_{\RR^2} U_{2,1}(x)\bigl(\cP(\rho_{2}) -
	\cP(\rho_{1})\bigr)(y)\dx\dy 
&\cr
&= - \int_{\RR^2}  \int_{\RR^2} U_{2,1}(x)
	\int_0^1 {{\rm d} \over \dlambda} \cP(\rho_\lambda)(y)\dlambda\dx\dy 
&\cr
& =  \beta \int_0^1 \left\langle \Bigl( U_{2,1} - 
\langle U_{2,1}\rangle(\lambda)\Bigr)^2\right\rangle (\lambda) \dlambda \, .
&\eqlbl\negdefness\cr}
$$
Since $\beta<0$, the last term in (\negdefness) is $\leq 0$, 
with $=0$ holding iff $U_{2,1}\equiv constant$. 
By (\negdefness) and (\posdefness) we conclude that $\rho_{2,1}\equiv 0$. 
Uniqueness is proved.\qed 

%\smallskip
%\noindent
%{\bf Remark 9.2:} {\it For certain $\Upsilon$, the 
%above uniqueness proof can be modified along the lines in 
%[\Aly, \KieLMP] to yield uniqueness for small $\beta >0$.}~\eqed

\smallskip
\noindent
{\bf Corollary 9.2:} {\it If $\beta <0$, $H \equiv const.$, and $\Upsilon$ 
is radially  symmetric, then the unique solution  $\rho_{\beta,H}$ of 
(\fixptEQ) is radially symmetric as well. 
Gauss' theorem then implies that the corresponding solution of (\PDE), 
$U_{H,\kappa}$ given in (\uofrho), is radial decreasing.}~\qed 
\smallskip

The proof of Corollary 9.2 is trivial. 
Theorem 9.1 and Corollary 9.2 prove Theorem~2.2.

% \vfill\eject

\bigskip
\noindent
{\bf ACKNOWLEDGMENT:} The work of Sagun Chanillo was supported
by NSF GRANT DMS-9623079, and the work of Michael Kiessling by
NSF GRANT DMS-9623220. 

\biblio

\bigskip
% \vfill\eject
\centerline{FIGURE CAPTIONS}
\bigskip
\bigskip
\bigskip
\noindent
Fig.1:\ Level curves  $e^{2u(x)}= 2^a$, $a\in\{-5,-4,...,0,1\}$,
for $u$ given by (\chakieSOL) with $n=2$, $|y|=1$, $\theta_0=0$, 
$\zeta =1$. $\max e^{2u}\approx 2.57$ is taken at the centers of 	
the two islands. For $|x|$ large, the conformal factor 
 $e^{2u(x)} \sim C |x|^{-8}$,  and the level curves become circular.
\bigskip
\bigskip
\noindent
Fig.2: Level curves $e^{2u}= 2^{a}$, $a\in\{-6,-5,...,0\}$, with $u$ 
given by (\stuartSOL), with $\zeta =1$, $y=- v^\prime$, $K_0=1$, 
$x_1=\langle x,v^\pr\rangle$ and $x_2=\langle x,v\rangle$.
$\max e^{2u}\approx 1.22$ is taken at the centers of the islands.
For $|\langle v,x\rangle|$ large, 
$e^{2u(x)} \sim C e^{-|\langle v,x\rangle|}$ and level curves become 
straight lines.

\bye